\documentclass[12pt,a4paper,reqno]{amsart}

\usepackage{enumerate}
\usepackage{float}
\usepackage{amssymb,amsmath}

\usepackage[all]{xy}
\newcommand{\Z}{{\mathbb Z}}

\newcommand{\F}{{\mathbb F}}

 \newcommand{\pcom}{_{p}^{\wedge}}

\renewcommand{\hom}{\operatorname{Hom}\nolimits}
\newcommand{\Mor}{\operatorname{Mor}\nolimits}
\newcommand{\Hom}{\operatorname{Hom}\nolimits}

\newcommand{\Aut}{\operatorname{Aut}\nolimits}%discrete groups
\newcommand{\Inn}{\operatorname{Inn}\nolimits}
\newcommand{\Out}{\operatorname{Out}\nolimits}

\newcommand{\End}{\operatorname{End}\nolimits}
\renewcommand{\ker}{\operatorname{Ker}\nolimits}

\newcommand{\Id}{\operatorname{Id}\nolimits}
\newcommand{\Syl}{\operatorname{Syl}\nolimits}
\newcommand{\rk}{\operatorname{rk}\nolimits}
\newcommand{\GL}{\operatorname{GL}\nolimits}
\newcommand{\SL}{\operatorname{SL}\nolimits}
\newcommand{\PGL}{\operatorname{PGL}\nolimits}
\newcommand{\PSL}{\operatorname{PSL}\nolimits}
\newcommand{\ASL}{\operatorname{ASL}\nolimits}

\newcommand{\res}{\operatorname{res}\nolimits}
\newcommand{\cor}{\operatorname{cor}\nolimits}
\newcommand{\Der}{\operatorname{Der}\nolimits}
\newcommand{\A}{\ifmmode{\mathcal{A}}\else${\mathcal{A}}$\fi}
\newcommand{\K}{\ifmmode{\mathcal{K}}\else${\mathcal{K}}$\fi}
\newcommand{\U}{\ifmmode{\mathcal{U}}\else${\mathcal{U}}$\fi}
\newcommand{\M}{\ifmmode{\mathcal{M}}\else${\mathcal{M}}$\fi}
\newcommand{\N}{\ifmmode{\mathcal{N}}\else${\mathcal{M}}$\fi}
\newcommand{\Ff}{\ifmmode{\mathcal{F}}\else${\mathcal{F}}$\fi}
\newcommand{\Ll}{\ifmmode{\mathcal{L}}\else${\mathcal{L}}$\fi}
\newtheorem{theorem}{Theorem}[section]
\newtheorem{proposition}[theorem]{Proposition}
\newtheorem{corollary}[theorem]{Corollary}
\newtheorem{lemma}[theorem]{Lemma}

\theoremstyle{definition}
\newtheorem{definition}[theorem]{Definition}
\newtheorem{remark}[theorem]{Remark}

\theoremstyle{remark}

\newtheorem*{Not}{Notation}
\newcommand{\definicio}{\stackrel{\text{def}}{=}}
\newcommand{\pb}{\parbox[h][6.5mm][c]{0mm}{}}
\newcommand{\pbb}{\parbox[h][8.5mm][c]{0mm}{}}
\newcommand{\pbbb}{\parbox[h][11mm][c]{0mm}{}}
\title[$p$-local finite groups of rank two for odd $p$]{All $p$-local finite groups of rank two for odd prime $p$}
\author{Antonio D{\'\i}az}
\address[A.~D{\'\i}az, A~Viruel]{
Departamento de {\'A}lgebra, Geometr{\'\i}a y Topolog{\'\i}a,
Universidad de M{\'a}\-la\-ga, Apdo correos 59, 29080 M{\'a}laga,
Spain.}
\email[A.~D{\'\i}az]{adiaz@agt.cie.uma.es} 
\email[A.~Viruel]{viruel@agt.cie.uma.es} 

\author{Albert Ruiz}
\address[A.~Ruiz]{Departament de Matem{\`a}tiques,
Universitat Aut{\`o}noma de Barcelona, 08193 Cerdanyola del
Vall{\`e}s, Spain.}
\email[A.~Ruiz]{Albert.Ruiz@uab.es}

\author{Antonio Viruel}
\thanks{
\indent First author is partially supported by MCED grant AP2001-2484. \\ 
\indent First and third authors are partially supported by grant MCYT-BFM2001-1825.\\
\indent Second author is partially supported by MCYT grant BFM2001-2035.}
\begin{document}

\begin{abstract}
In this paper we give a classification of the rank two $p$-local finite groups
for odd~$p$. This study requires the analisis of the possible saturated
fusion systems in terms of the outer automorphism group ant the proper 
$\Ff$-radical subgroups.
Also, for each case in the classification, either we give a finite group with the corresponding
fusion system or we check that it corresponds to an exotic $p$-local
finite group, getting some new examples of these for~$p=3$.
\end{abstract}

\maketitle

%%%%%%%%%%%%%%%%%%%%%%%%%%%%%%%%%%%%%%%%%%%%%
% Introduction
%%%%%%%%%%%%%%%%%%%%%%%%%%%%%%%%%%%%%%%%%%%%%

%\input{section1}
\section{Introduction}

When studying the $p$-local homotopy theory of classifying spaces
of finite groups, Broto-Levi-Oliver \cite{BLO2} introduced the
concept of $p$-local finite group as a $p$-local analogue of the
classical concept of finite group. These purely algebraic objects,
whose basic properties are reviewed in Section
\ref{p-local-review}, are a generalization of the classical theory
of finite groups in the sense that every finite group leads to a
$p$-local finite group, although there exist exotic $p$-local
finite groups which are not associated to any finite group as it
can be read in \cite[Sect.\ 9]{BLO2}, \cite{LO}, \cite{RV} or
Theorem \ref{sfs-B-gamma} below. Besides its own interest, the
systematic study of possible $p$-local finite groups, i.e.\
possible $p$-local structures, is meaningful when working in other
research areas as transformation groups (e.g.\ when constructing
actions on spheres \cite{adem-smith,adem-fixity}), or modular
representation theory (e.g.\ the study of the $p$-local structure
of a group is a very first step when verifying conjectures like
those of Alperin \cite{Alp-conj} or Dade \cite{Dade-conj}). It
also provides an opportunity to enlighten one of the highest
mathematical achievements in the last decades: The Classification
of Finite Simple Groups \cite{gorenstein-lyons-solomon}. A
milestone in the proof of that classification is the
characterization of finite simple groups of $2$-rank two (e.g.\
\cite[Chapter 1]{gorenstein-83}) what is based in a deep
understanding of the $2$-fusion of finite simple groups of low
$2$-rank \cite[1.35, 4.88]{gorenstein-82}. Unfortunately, almost
nothing seems to be known about the $p$-fusion of finite groups of
$p$-rank two for an odd prime $p$ \cite{dietz-priddy}, and this
work intends to remedy that lack of information by classifying all
possible saturated fusion systems over finite $p$-groups of rank
two, $p>2$.

\begin{theorem}\label{teo:elteorema}
Let $p$ be an odd prime and $S$ be a rank two $p$-group. Given a saturated fusion system
$(S,\Ff)$, one of the following hold:
\begin{itemize}
\item $\Ff$ has no proper $\Ff$-centric, $\Ff$-radical subgroups
and it corresponds to the group $S:\Out_\Ff(S)$.
\item $\Ff$ is the fusion system of a group $G$ which fits in the following
extension:
$$
1 \to S_0 \to G \to W \to 1
$$
where $S_0$ is a subgroup of index $p$ in $S$ and $W$ contains $\SL_2(p)$.
\item $\Ff$ is the fusion system of an extension of one of the following
finite simple groups:
  \begin{itemize}
  \item $L_3(p)$ for any $p$,
  \item ${}^2F_4(2)'$, $J_4$, $L_3^{\pm}(q)$, ${}^3D_4(q)$, ${}^2F_4(q)$ for $p=3$,
  where $q$ is a $3'$ prime power,
  \item $Th$ for $p=3$,
  \item $He$, $Fi'_{24}$, $O'N$ for $p=7$ or
  \item $M$ for $p=13$.
  \end{itemize}
\item $\Ff$ is an exotic fusion system characterized by the following data:
  \begin{itemize}
  \item $S=7^{1+2}_+$ and all the rank two elementary abelian subgroups are $\Ff$-radical,
  \item $S=G(3,2k;0,\gamma,0)$ and the only proper $\Ff$-radical, $\Ff$-centric subgroups
  are one or two $S$-conjugacy classes of rank two elementary abelian subgroups,
  \item $S=G(3,2k+1;0,0,0)$ and the proper $\Ff$-centric, $\Ff$-radical subgroups are
  either one or two $S$-conjugacy classes of subgroups isomorphic to $3^{1+2}_+$,
  either one $S$-conjugacy class of rank two elementary abelian subgroups and one
  subgroup isomorphic to $\Z/3^k+\Z/3^k$.
  \end{itemize}
\end{itemize}
\end{theorem}

\begin{proof}
For an odd prime $p$, the isomorphism type of a rank two
$p$-group, namely $S$, is described in Theorem
\ref{thm:clasificacionrango2}, hence the proof is done by studying
the different cases of $S$.

For $p>3$, Theorems \ref{M-resistant}, \ref{C-resistant} and
\ref{Thm:g(r,e)isresistant} show that any saturated fusion system
$(S,\Ff)$ is induced by $S:\Out_\Ff(S)$, i.e.\ $S$ is resistant
(see Definition \ref{def:resistant}), unless $S\cong p^{1+2}_+$,
hence \cite{RV} completes the proof for the $p>3$ case.

For $p=3$, Theorems \ref{M-resistant}, \ref{C-resistant},
\ref{B(n;0,b,c,d)resistant} and \ref{Thm:B(r;0,1,0,0)resistant}
describe all rank two $3$-groups which are resistant. The
saturated fusion system over non resistant rank two $3$-groups are
then obtained in Theorems \ref{sfs-G-3} and \ref{sfs-B-gamma} when
$S\not\cong 3^{1+2}_+$, what completes the information in
\cite{RV} and finishes the proof.
\end{proof}

It is worth noticing that along the proof of the theorem above
some interesting contributions are made:
\begin{itemize}
\item The Appendix provides a neat compendium of the group theoretical
properties of rank two $p$-groups, $p$ odd, including a
description of their automorphism groups, and $p$-centric
subgroups. It does not only collect the related results in the
literature \cite{BB,BB1,dietz,dietz-priddy,huppert}, but
extends them.

\item  Theorems \ref{M-resistant}, \ref{C-resistant},
\ref{Thm:g(r,e)isresistant}, \ref{B(n;0,b,c,d)resistant} and
\ref{Thm:B(r;0,1,0,0)resistant} identify a large family of
resistant groups complementing the results in \cite{stancu}, and
extending those in \cite{martino-priddy}.

\item Theorem \ref{sfs-B-gamma} provides infinite families
of exotic rank two $p$-local finite groups with arbitrary large
Sylow $p$-subgroup. Unlikely the other known infinite families of exotic
$p$-local finite groups \cite{BM,LO}, some of these new families
cannot be constructed as the homotopy fixed points of automorphism
of a $p$-compact group. Nevertheless, it is still possible to
construct an ``ascending" chain of exotic $p$-local finite groups
whose colimit, we conjecture, should provide an example of an
exotic $p$-local compact group \cite[Section 6]{BLOs}.
\end{itemize}

\noindent{\bf Organization of the paper:} Along Section
\ref{p-local-review} we quickly review the basics on $p$-local
finite groups. In Section \ref{resistant} we define the concept of
resistant $p$-group, similar to that of Swan group, and develop
some machinery to identify resistant groups. In Section
\ref{metacyclic} we study the fusion systems over non maximal
nilpotency class rank two $p$-groups, while the study of fusion
systems over maximal nilpotency class rank two $p$-groups is done
in Section \ref{maximal_class}. We finish this paper with an
Appendix collecting the group theoretical background on rank two
$p$-groups which is needed along the classification.

\noindent{\bf Notation:} By $p$ we always denote an odd prime, and
$S$ a $p$-group of size $p^r$. The group theoretical notation used
along this paper is that described in the Atlas \cite[5.2]{atlas}.
For a group~$G$, and $g\in G$, we denote by $c_g$ the conjugation
morphism $x\mapsto g^{-1}xg$. If $P,Q\leq G$, the set of
$G$-conjugation morphisms from $P$ to $Q$ is denoted by
$\Hom_G(P,Q)$, so if $P=Q$ then $\Aut_G(P)=\Hom_G(P,P)$. Notice
that $\Aut_G(G)$ is then $\Inn(G)$, the group of inner
automorphisms of $G$. Given $P$ and $Q$ groups,
$\operatorname{Inj(P,Q)}$ denotes the set of injective
homomorphisms.

\noindent{\bf Acknowledges:} The authors are indebted
with Ian Leary, Avinoam Mann and Bob Oliver. The first two
provided us the references for stating Theorem
\ref{thm:clasificacionrango2}, while the later has shown
interest and made helpful comments and suggestions throughout the
course of this work.

%%%%%%%%%%%%%%%%%%%%%%%%%%%%%%%%%%%%%%%%%%%%%
% 2. p-local finite groups
%%%%%%%%%%%%%%%%%%%%%%%%%%%%%%%%%%%%%%%%%%%%%

%\input{section2}
\section{$p$-local finite groups}\label{p-local-review}

At the beginning of this section we quickly review the concept of
$p$-local finite group introduced in \cite{BLO2} that builds on a
previous unpublished work of L.~Puig, where the axioms for fusion
systems are already established. See \cite{BLOs} for a survey on
this subject.

After that we study some particular cases where a saturated
fusion system is controlled by the normalizer of the Sylow
$p$-subgroup.

We end this section with some tools which allow us to study the
exoticism of an abstract saturated fusion system.

\begin{definition}\label{def:fusion_system}
A fusion system $\Ff$ over a finite $p$-group $S$ is a category
whose objects are the subgroups of $S$, and whose morphisms sets
$\Hom_\Ff(P,Q)$ satisfy the following two conditions:
\begin{enumerate}[(a)]
\item $\Hom_S(P,Q) \subseteq \Hom_\Ff(P,Q) \subseteq
\operatorname{Inj}(P,Q)$ for all $P$ and $Q$ subgroups of $S$.

\item Every morphism in $\Ff$ factors as an isomorphism in $\Ff$
followed by an inclusion.
\end{enumerate}
\end{definition}

We say that two subgroups $P$,$Q \leq S$ are
\textit{$\Ff$-conjugate} if there is an isomorphism between them
in $\Ff$. As all the morphisms are injective by condition (b), we
denote by $\Aut_\Ff(P)$ the group $\Hom_\Ff(P,P)$. We denote by
$\Out_\Ff(P)$ the quotient group $\Aut_\Ff(P)/\Aut_P(P)$.

The fusion systems that we consider are saturated, so we need the following
definitions:

\begin{definition}\label{def:N}
Let $\Ff$ be a fusion system over a $p$-group $S$.
\begin{itemize}
\item A subgroup $P \leq S$ is \emph{fully centralized in $\Ff$}
if $|C_S(P)|\geq|C_S(P')|$ for all $P'$ which is $\Ff$-conjugate
to $P$. \item A subgroup $P \leq S$ is \emph{fully normalized in
$\Ff$} if $|N_S(P)|\geq|N_S(P')|$ for all $P'$ which is
$\Ff$-conjugate to $P$. \item $\Ff$ is a \emph{saturated fusion
system} if the following two conditions hold:
\begin{enumerate}[(I)]
\item Every fully normalized subgroup $P \leq S$ is fully
centralized and $\Aut_S(P) \in \Syl_p(\Aut_\Ff(P))$. \item If
$P\leq S$ and $\varphi \in \Hom_\Ff(P,S)$ are such that $\varphi
P$ is fully centralized, and if we set
$$ N_\varphi=\{ g \in N_S(P) \mid \varphi c_g \varphi^{-1} \in
\Aut_S(\varphi P)\} ,
$$ then there is $\overline\varphi \in \Hom_\Ff(N_\varphi,S)$ such
that $\overline\varphi|_P=\varphi$.
\end{enumerate}
\end{itemize}
\end{definition}

\begin{remark}\label{rmk:orderW}
From the definition of fully normalized and the condition (I) of
saturated fusion system we get that if $\Ff$ is a saturated fusion
system over a $p$-group $S$ then $p$ cannot divide the order of the outer
automorphism group $\Out_\Ff(S)$.
\end{remark}

As expected, every finite group $G$ gives rise to a saturated
fusion system \cite[Proposition 1.3]{BLO2}, which provides
valuable information about  $BG\pcom$ \cite{O1}. Some classical
results for finite groups can be generalized to saturated fusion
systems, as for example, Alperin's fusion theorem for saturated fusion systems \cite[Theorem
A.10]{BLO2}:

\begin{definition}\label{def:centric-radical}
Let $\Ff$ be any fusion system over a $p$-group $S$. A subgroup
$P\leq S$ is:
\begin{itemize}
\item  \emph{$\Ff$-centric} if $P$ and all its $\Ff$-conjugates
contain their $S$-centralizers.
\item \emph{$\Ff$-radical} if $\Out_\Ff(P)$ is
$p$-reduced, that is, if $\Out_\Ff(P)$ has no nontrivial normal $p$-subgroups.
\end{itemize}
\end{definition}

\begin{theorem}[Alperin's fusion theorem for saturated fusion systems]\label{teo:alperin}
Let $\Ff$ be a saturated fusion system over $S$. Then for each
morphism $\psi\in\Aut_\Ff(P,P')$, there exists a sequence of
subgroups of $S$
$$P=P_0, P_1,\ldots,P_k=P'\quad\text{and}\quad Q_1,Q_2,\ldots,Q_k,$$
and morphisms $\psi_i\in\Aut_\Ff(Q_i)$, such that
\begin{itemize}
\item $Q_i$ is fully normalized in $\Ff$, $\Ff$-radical and
  $\Ff$-centric for each $i$;
\item $P_{i-1},P_i\leq Q_i$ and $\psi_i(P_{i-1})=P_i$ for each
$i$;
  and
\item $\psi=\psi_k\circ\psi_{k-1}\circ\cdots\circ\psi_1$.
\end{itemize}
\end{theorem}
%This theorem points up that the $\Ff$-centric, $\Ff$-radical subgroups
%will play an important role on the classification, so, for shortness, we
%will use the next definition:

The subgroups $Q_i$'s in the theorem above determine the structure of $\Ff$, so they deserve a
name:

\begin{definition}\label{def:F-Alperin}
Let $\Ff$ be any fusion system over a $p$-group $S$. We say that a
subgroup~$Q \leq S$ is \emph{$\Ff$-Alperin} if it is fully
normalized in $\Ff$, $\Ff$-radical and $\Ff$-centric.
\end{definition}

The definition of $p$-local finite group still requires one more
new concept.

Let $\Ff^c$ denote the full
subcategory of $\Ff$  whose objects are the $\Ff$-centric
subgroups of $S$.

\begin{definition}\label{def:centric_linking_system}
Let $\Ff$ be a fusion system over the $p$-group $S$. A
\emph{centric linking  system associated to $\Ff$} is a category
$\Ll$ whose objects are the $\Ff$-centric subgroups of $S$,
together with a functor
$$\pi\colon\Ll\longrightarrow\Ff^c$$ and ``distinguished''
monomorphisms $P{\buildrel{\delta_P}\over{\longrightarrow}}
\Aut_\Ll(P)$ for each $\Ff$-centric subgroup $P\leq S$, which
satisfy the following conditions:
\begin{itemize}
\item[{\rm (A)}] $\pi$ is the identity on objects and surjective
on morphisms. More precisely, for each pair of objects
$P,Q\!\in\!\Ll$, $Z(P)$ acts freely on $\Mor_\Ll(P,Q)$ by composition
(upon identifying $Z(P)$ with $\delta_P(Z(P))\leq\Aut_\Ll(P)$),
and $\pi$ induces a bijection
$$\Mor_\Ll(P,Q)/Z(P){\buildrel{\cong}\over{\longrightarrow}}\Hom_\Ff(P,Q).$$

\item[{\rm (B)}] For each $\Ff$-centric subgroup $P\leq S$ and
each
  $g\in P$, $\pi$ sends $\delta_P(g)\in\Aut_\Ll(P)$ to
  $c_g\in\Aut_\Ff(P)$.

\item[{\rm (C)}] For each $f\in\Mor_\Ll(P,Q)$ and each $g\in P$,
the
  equality $\delta_Q(\pi(f)(g))\circ f=f\circ\delta_P(g)$ holds in
  $\Ll$
\end{itemize}
\end{definition}

Finally, the definition of $p$-local finite group is:

\begin{definition}\label{p-local-def}
A \emph{$p$-local finite group} is a triple $(S,\Ff,\Ll)$, where
$S$ is a $p$-group, $\Ff$ is a saturated fusion system over $S$
and $\Ll$ is a centric linking system associated to $\Ff$. The
\emph{classifying space} of the $p$-local finite group
$(S,\Ff,\Ll)$ is the space $|\Ll|\pcom$.
\end{definition}

Given a fusion system $\Ff$ over the $p$-group $S$, there exists
an obstruction theory for the existence and uniqueness of a
centric linking system, i.e. a $p$-local finite group, associated
to~$\Ff$. The question is solved for $p$-groups of small rank by
the following result \cite[Theorem~E]{BLO2}:

\begin{theorem}
Let $\Ff$ be any saturated fusion system over a $p$-group $S$.  If
$\rk_p(S)<p^3$, then there exists a centric linking system
associated to $\Ff$. And if $\rk_p(S)<p^2$, then there exists a
unique centric linking system associated to $\Ff$.
\end{theorem}

As all $p$-local finite groups studied in this work are over rank
two $p$-groups $S$, we obtain:

\begin{corollary}\label{cor:rank2}
Let $p$ be an odd prime. Then the set of $p$-local finite groups
over a rank two $p$-group $S$ is in bijective correspondence with
the set of saturated fusion systems over $S$.
\end{corollary}

In \cite[Section 2]{BLO2} is defined the ``centralizer'' fusion system of a
given fully centralized subgroup:
\begin{definition}\label{def:centralizer}
Let $\Ff$ be a fusion system over $S$ and $P \leq S$ a fully
centralized subgroup in $\Ff$. The \emph{centralizer fusion system
of $P$ in $\Ff$, $C_\Ff(P)$} is the fusion system over $C_S(P)$
with objects $Q\leq C_S(P)$ and morphisms, $$
\begin{array}{l}
\Hom_{C_\Ff(P)}(Q,Q') \\
\quad =\{ \varphi  \in \Hom_{\Ff}(Q,Q') | \exists
\psi \!\in\! \Hom_{\Ff}(QP,Q'P), \psi |_Q=\varphi ,\psi |_P =\Id_P
\} \,.
\end{array}
$$
\end{definition}

\begin{remark}\label{rmk:centralizadorygrupo}
If we consider the fusion system corresponding to a finite group
$G$ with Sylow $p$-subgroup $S$ and $P \leq S$ such that
$C_S(P)\in\Syl_p(C_G(P))$, i.e.\ $P$ is fully centralized in
$\Ff_S(G)$, then the fusion system $(C_S(P),C_\Ff(P))$ is the
fusion system $(C_S(P),\Ff_{C_S(P)}(C_G(P)))$, so it is again the
fusion system of a finite group.
\end{remark}

In \cite[Section 6]{BLO2} the ``normalizer'' fusion system of a
given fully normalized subgroup is defined:

\begin{definition}
Let $\Ff$ be a fusion system over $S$ and $P\leq S$ a fully
normalized subgroup in $\Ff$. The \emph{normalizer fusion system
of $P$ in $\Ff$, $N_\Ff(P)$} is the fusion system over $N_S(P)$
with objects $Q\leq N_S(P)$ and morphisms, $$
\begin{array}{l}
\Hom_{N_\Ff(P)}(Q,Q') \\
 \quad =\{ \varphi  \in \Hom_{\Ff}(Q,Q') | \exists
\psi \!\in\! \Hom_{\Ff}(QP,Q'P), \psi |_Q=\varphi ,\psi |_P \!\in\! \Aut(P)
\} \,.
\end{array}
$$
\end{definition}

Also in \cite[Section 6]{BLO2} it is proved that if $\Ff$ is a
saturated fusion system over $S$ and $P$ is a fully normalized
subgroup then $N_\Ff(P)$ is a saturated fusion system over
$N_S(P)$.

\begin{remark}\label{rmk:sylow-controla-por-def}
For a fusion system $\Ff$ over $S$, when considering the
normalizer fusion system over the own Sylow $S$, it turns out that
$N_\Ff(S)=\Ff$ if and only if every $\varphi \in \Hom_\Ff(Q,Q')$
extends to $\psi \in \Aut_\Ff(S)$ for each $Q,Q'\leq S$.
\end{remark}

Moreover, we have the following two characterizations of $\Ff$
reducing to the normalizer of the Sylow:
\begin{lemma}\label{thm:sylow-controla-1}
Let $\Ff$ be a saturated fusion system over the $p$-group $S$.
Then $\Ff=N_\Ff(S)$ if and only if $S$ itself is the only
$\Ff$-Alperin subgroup of $S$.
\end{lemma}
\begin{proof}
If $S$ is the unique $\Ff$-Alperin subgroup then
the assertion follows from Alperin's fusion theorem for saturated
fusion systems (Theorem \ref{teo:alperin}). Assume then that
$\Ff=N_\Ff(S)$ and choose $P\leq S$ $\Ff$-Alperin. Using that
$\Ff=N_\Ff(S)$, it is straightforward that $\Out_S(P)$ is normal
in $\Out_\Ff(P)$ and, as $P$ is $\Ff$-radical, it must be trivial.
Then $P=N_S(P)$, and as $S$ is a $p$-group, $P$ must be equal to
$S$.
\end{proof}

\begin{lemma}\label{thm:sylow-controla-2}
Let $\Ff$ be a saturated fusion system over the $p$-group $S$. Then
$\Ff=N_\Ff(S)$ if and only if $N_\Ff(P)=N_{N_\Ff(P)}(N_S(P))$ for
every $P\leq S$ fully normalized in $\Ff$.
\end{lemma}
\begin{proof}
Assume first that $\Ff=N_\Ff(S)$ and $P\leq S$ is fully normalized
in $\Ff$. In general we have that $N_\Ff(P) \supset
N_{N_\Ff(P)}(N_S(P))$. The other inclusion follows if all $\varphi
\in \Hom_{N_\Ff(P)}(Q,Q')$ which extend to a morphism $\psi \in
\Hom_\Ff(PQ,PQ')$ such that $\psi |_Q=\varphi$ and $\psi |_P \in
\Aut_\Ff(P)$, extend to an element of $\Aut(N_S(P))$. But using
Remark \ref{rmk:sylow-controla-por-def} we get that $\psi$ extends
to an element of $\Aut_\Ff(S)$, which restricts to an element of
$\Aut_\Ff(N_S(P))$ because $\psi$ restricts to an element of
$\Aut(P)$.

Assume now that for every $P\leq S$ fully normalized in $\Ff$ we
have $N_\Ff(P)=N_{N_\Ff(P)}(N_S(P))$. According to Lemma
\ref{thm:sylow-controla-1} we have to check that $S$ does not
contain any proper $\Ff$-Alperin subgroup. Let $P$ be a
$\Ff$-Alperin subgroup, then it is $N_\Ff(P)$-Alperin too. But,
applying Lemma \ref{thm:sylow-controla-1} to $N_\Ff(P)$, we get
that $S=P$.
\end{proof}

Finally in this section we give some results which allow us to
determine in some special cases the existence of a finite group
with a fixed saturated fusion system.

We begin with a definition which does only depend on the $p$-group $S$:

\begin{definition}
Let $S$ be a $p$-group. A subgroup $P\leq S$ is
\emph{$p$-centric} in $S$ if $C_S(P)=Z(P)$.
\end{definition}

\begin{remark}\label{rmk:p-centric}
If $\Ff$ is any fusion system over the $p$-group $S$, then
$\Ff$-centric subgroups are $p$-centric subgroups too.
\end{remark}

The following result is a generalization of \cite[Proposition
9.2]{BLO2} which apply to some of our cases. Recall that given a
fusion system $(S,\Ff)$, a subgroup $P\leq S$ is called
\emph{strongly closed in $\Ff$} if no element of $P$ if
$\Ff$-conjugate to any element of $S\setminus P$.

\begin{proposition}\label{proposition9.2deBLO2}
Let $(S,\Ff)$ be a saturated fusion system such that every non
trivial strongly closed subgroup $P\leq S$ is non elementary
abelian, $p$-centric and does not factorize as a product of two or
more subgroups which are permuted transitively by $\Aut_\Ff(P)$.
Then if $\Ff$ is the fusion system of a finite group, it is the
fusion system of a finite almost simple group.
\end{proposition}
\begin{proof}
Suppose $\Ff=\Ff_S(G)$ for a finite group $G$. Assume also that
$\#|G|$ is minimal with this property. Consider $H \lhd G$ a
minimal non trivial normal subgroup in $G$. Then $H\cap S$ is a
strongly closed subgroup of $(S,\Ff)$. If $H\cap S=1$ then $\Ff$
is also the fusion system of $G/H$, which contradicts the
assumption of minimality on $G$. Now $P\definicio H\cap S$ is a
non trivial normal strongly closed subgroup in $(S,\Ff)$ and, as
$H$ is normal, $P$ is the Sylow $p$-subgroup in $H$. By
\cite[Theorem 2.1.5]{gorenstein-68} $H$ must be either elementary
abelian or a product of non abelian simple groups isomorphic one
to each other which must be permuted transitively by $N_G(H)=G$
(now by minimality of $H$). Notice that $H$ cannot be elementary
abelian as that would imply that $P$ is so while this is not
possible by hypothesis. Therefore $H$ is a product of non abelian
simple groups. If $H$ is not simple (so there is more than one
factor) this would break $P$ into two ore more factors which would
be permuted transitively, so $H$ must be simple. Now, as $P$ is
$p$-centric, $C_G(H)\cap S \subset P \subset H$, so $C_G(H)\cap S
\subset C_G(H)\cap H=1$ and $C_G(H)=1$ ($C_G(H)$ is a normal
subgroup in $G$ of order prime to $p$, so if $C_G(H)\neq 1$ taking
$G/C_G(H)$ gives again a contradiction with the minimality of
$\#|G|$). This tells us that $H \lhd G \leq \Aut(H)$, so $G$ is
almost simple.
\end{proof}

\begin{remark}
In fact \cite[Proposition 9.2]{BLO2} proves that if the only non
trivial strongly closed $p$-subgroup in $(S,\Ff)$ is $S$, and
moreover $S$ is non abelian and it does not factorize as a product
of two or more subgroups which are permuted transitively by
$\Aut_\Ff(S)$, then if $(S,\Ff)$ is the fusion system of a finite
group, it is the fusion system of an extension of a simple group
by outer automorphisms of order prime to $p$.
\end{remark}

We finish this section with the following result, which can be found in
\cite[Corollary 6.17]{BCGLO2}:

\begin{lemma}\label{lem:extensionescentrales}
Let $(S,\Ff)$ be a saturated fusion system, and assume there
is a nontrivial subgroup $A\leq Z(S)$ which is central in $\Ff$
(i.e. $C_\Ff(A)=\Ff$).
Then $\Ff$ is the fusion system of a finite group if and only if
$\Ff/A$ is so.
\end{lemma}

%%%%%%%%%%%%%%%%%%%%%%%%%%%%%%%%%%%%%%%
% 3. Resistant p-groups
%%%%%%%%%%%%%%%%%%%%%%%%%%%%%%%%%%%%%%%

%\input{section3}
\section{Resistant $p$-groups}\label{resistant}

In this section we recall the notion of Swan group and introduce
its generalization for fusion systems, as well as we discuss some
related results. Some of these results were considered
independently by Stancu in \cite{StancuThesis}.

For a fixed prime $p$ we recall that a subgroup $H\leq G$ is said
to \emph{control (strong) $p$-fusion} in~$G$ if $H$ contains a
Sylow $p$-subgroup of $G$ and whenever $P$, $g^{-1}Pg\leq H$ for $P$
a $p$-subgroup of $G$ then $g=hc$ where $h\in H$ and $c\in
C_G(P)$.

If we focus on a $p$-group $S$ we may wonder if \emph{$N_G(S)$
controls $p$-fusion} in $G$ whenever $S\in \Syl_p(G)$. If this is
always the case, then $S$ is called  a \emph{Swan} group. The
equivalent concept in the setting of $p$-local finite groups is:

\begin{definition}\label{def:resistant}
A $p$-group $S$ is called \emph{resistant} if whenever $\Ff$ is a
saturated fusion system over $S$ then $N_\Ff(S)=\Ff$.
\end{definition}

\begin{remark}\label{rmk:resistant-swan}
Considering the saturated fusion system associated to $S\in
\Syl_p(G)$ it is clear that every resistant group is a Swan group.
In the opposite way, up to date there is no known Swan group that
is not a resistant group.
\end{remark}

Following \cite[Section 2]{martino-priddy} we look for conditions
on a $p$-group $S$ for being resistant. If the $p$-group $S$ is
resistant, when treating with a fusion system $\Ff$ over $S$ all
morphisms in $\Ff$ are restrictions of $\Ff$-automorphisms of $S$.
In the general case, we must pay attention to possible
$\Ff$-Alperin subgroups to understand the whole category $\Ff$.
The first step towards this objective is to examine
$p$-centric subgroups.

\begin{theorem}\label{thm:suficiente-para-resistant}
Let $S$ be a $p$-group. If every proper $p$-centric subgroup $P\lneq S$
verifies $$ \Out_S(P)\cap O_p(\Out(P))\neq
1 \, , $$ then $S$ is a resistant group.
\end{theorem}
\begin{proof}
Let $\Ff$ be a saturated fusion system over $S$. According to
Lemma \ref{thm:sylow-controla-1}, it is enough to prove that $S$
is the only $\Ff$-Alperin subgroup. Let $P\lneq S$ be a proper
$\Ff$-Alperin subgroup, hence $p$-centric by Remark
\ref{rmk:p-centric}. As $\Out_S(P)\leq \Out_\Ff(P)$, we have that
$$ 1\neq \Out_S(P)\cap O_p(\Out(P))\leq \Out_\Ff(P)\cap
O_p(\Out(P)) \trianglelefteq O_p(\Out_\Ff(P))\, . $$ Hence $P$ cannot be $\Ff$-radical.
\end{proof}

We obtain a family of resistant groups:

\begin{corollary}\label{cor:abelian}
Abelian $p$-groups are resistant groups.
\end{corollary}

Also is meaningful determining whether or not a $p$-group can be
$\Ff$-Alperin for some saturated fusion system $\Ff$:

\begin{lemma}\label{lem:suficiente-para-no-radical}
Let $P$ be a $p$-group such that $O_p(\Out(P))$ is the Sylow
$p$-subgroup of $\Out(P)$. Then $P$ is not $\Ff$-Alperin
for any saturated fusion system $\Ff$ over $S$ with $P\lneq S$.
\end{lemma}
\begin{proof}
Let $S$ be a $p$-group with $P\lneq S$ and $\Ff$ be a
saturated fusion system over $S$. If $P$ is $\Ff$-radical then the
normal $p$-subgroup $O_p(\Out(P))\cap \Out_\Ff(P)$ of
$\Out_\Ff(P)$ must be trivial. Being $O_p(\Out(P))$ a Sylow
$p$-subgroup and normal this implies that $\Out_\Ff(P)$ is a
$p'$-group and so, by definition, $\Aut_P(P)\in
\Syl_p(\Aut_\Ff(P))$. If in addition $P$ is fully normalized then
we know that $\Aut_S(P)$ is another Sylow $p$-subgroup of
$\Aut_\Ff(P)$. So they both must be same size. Finally, if $P$ is
$\Ff$-centric then $Z(P)=C_S(P)$ and $P=N_S(P)$, which is false
for $p$-groups unless $P$ is equal to $S$.
\end{proof}

We obtain some useful corollaries:

\begin{corollary}\label{cor:3.7}
$\Z/p^n$ is not $\Ff$-Alperin
in any saturated fusion system $\Ff$ over $S$ where $\Z/p^n\lneq
S$.
\end{corollary}
\begin{proof}
An easy verification shows that $\Aut(\Z/p^n)$ is equal to
$(\Z/p^n)^*$, which is abelian. Now just apply Lemma
\ref{lem:suficiente-para-no-radical}.
\end{proof}

\begin{corollary}\label{cor:3.6} If $n$ and $m$ are two different non zero
positive integers, then
$\Z/p^n\times \Z/p^m$ is not
$\Ff$-Alperin in any saturated fusion system $\Ff$ over $S$ where
$\Z/p^n\times \Z/p^m\lneq S$.
\end{corollary}
\begin{proof}
Suppose $n<m$ and take $f\in \End(\Z/p^n\times \Z/p^m)$ with
$f(\bar{1},\bar{0})=(\bar{a},\bar{b})$ and
$f(\bar{0},\bar{1})=(\bar{c},\bar{d})$. It must hold that $b \equiv 0$
mod $p^{m-n}$. Moreover, $f$ is an automorphism if and only if
order$(\bar{a},\bar{b})=p^n$, order$(\bar{c},\bar{d})=p^m$ and
$\langle(\bar{a},\bar{b}),(\bar{c},\bar{d})\rangle=\Z/p^n\times \Z/p^m$. It
can be checked that the first condition is equivalent to $a\neq 0$
mod $p$, the second to $d\neq 0$ mod $p$, and the third is
consequence of the previous ones. Counting elements it turns out
that $\Aut(\Z/p^n\times \Z/p^m)$ has order $p^{3n+m-2}(p-1)^2$.
Finally, the subgroup $\{(\bar{a},\bar{b}),(\bar{c},\bar{d})$ with
$a=1$ mod $p$ and $d=1$ mod $p\}$ is normal and has order
$p^{3n+m-2}$.
\end{proof}

\begin{corollary}\label{cor:metanonabelian}
If $p$ is odd then non abelian metacyclic $p$-groups $M$ are not
$\Ff$-Alperin in any saturated fusion system $\Ff$ over $S$ with
$M\lneq S$.
\end{corollary}
\begin{proof}
According to \cite[Section 3]{dietz}, for $p$ odd and $M$ a non
abelian metacyclic $p$-group, $O_p(\Out(M))$ is the Sylow
$p$-subgroup of $\Out(M)$, so the result follows using Lemma
\ref{lem:suficiente-para-no-radical}.
\end{proof}

Recall the notation in Theorem \ref{thm:clasificacionrango2} for the families
of $p$-rank two $p$-groups and Theorem~\ref{thm:maxclass} for the maximal
nilpotency class $3$-rank two $3$-groups.

\begin{corollary}\label{cor:Grenopuede}
If $p$ is odd then $G(p,r;\epsilon)$ cannot be $\Ff$-Alperin in any saturated
fusion system over $S$, with $G(p,r;\epsilon) \lneq S$.
\end{corollary}
\begin{proof}
By \cite[Proposition 1.6]{dietz-priddy}  $O_p(\Out(G(p,r;\epsilon)\!)\!)$ is the
Sylow $p$-subgroup of 
$\Out(G(p,r;\epsilon))$, so applying Lemma \ref{lem:suficiente-para-no-radical}
$G(p,r;\epsilon)$ could not be $\Ff$-Alperin.
\end{proof}

\begin{corollary}\label{B(n)_noes}
$B(3,r;\beta,\gamma,\delta)$ is not $\Ff$-Alperin in any saturated fusion system $\Ff$
over $S$ where $B(3,r;\beta,\gamma,\delta)\lneq S$.
\end{corollary}
\begin{proof}
From \cite{BB} we have that the Frattini subgroup of $B(3,r;\beta,\gamma,\delta)$,
$\Phi(B(3,r;\beta,\gamma,\delta))$, is $\langle s_2,s_3,\dots,s_{r-1} \rangle = \langle s_2,s_3
\rangle$. Consider the Frattini map
$$
B(3,r;\beta,\gamma,\delta)) \rightarrow
B(3,r;\beta,\gamma,\delta)/\Phi(B(3,r;\beta,\gamma,\delta))\simeq \langle \overline{s},\overline{s_1} \rangle
$$
and the induced map
$$
\rho:\Out(G) \rightarrow \Aut(B(3,r;\beta,\gamma,\delta)/\Phi(B(3,r;\beta,\gamma,\delta)))\simeq \GL_2(3)
$$
whose kernel is a $3$-group. As $\gamma_1$ is characteristic in
$B(3,r;\beta,\gamma,\delta)$, for every class of morphisms~$\varphi$ in $\Out(G)$ it
holds that $\rho(\varphi)(\langle\overline{s_1}\rangle)\leq
\langle\overline{s_1}\rangle$, and so the image of $\rho$ is
contained in the lower triangular matrices.

The subgroup generated by $\big(\!\begin{smallmatrix} 1 & 0
\\ 1 & 1\end{smallmatrix}\!\big)$ is normal and is the Sylow
$3$-subgroup of the lower triangular matrices of $\GL_2(3)$. Its
preimage by $\rho$ is normal in $\Out(G)$ and, as the kernel of
$\rho$ is a $3$-group, it is the Sylow $3$-subgroup of $\Out(G)$
too. To finish the proof apply Lemma
\ref{lem:suficiente-para-no-radical}.
\end{proof}

The following result related to Corollary \ref{cor:3.6} is needed
in successive sections, but before we recall here \cite[Lemma
4.1]{RV} in a clearer form:

\begin{lemma}\label{util}
Let $G$ be a $p$-reduced subgroup (that is, $G$ has no nontrivial
normal $p$-subgroup) of $\GL_2(p)$, $p\geq3$. If $p$ divides the
order of $G$ then $\SL_2(p)\leq G$.
\end{lemma}
\begin{proof}
If $p$ divides the order of the group $G$, then there is an
element of order $p$. As $G$ is $p$-reduced it cannot be the only
one, so we have a subgroup of $\GL_2(p)$ with more than one
nontrivial $p$-subgroup. Observe that the only nontrivial
$p$-subroups in $\GL_2(p)$ are the Sylow $p$-subgroups, so $G$ has
more than one Sylow $p$-subgroup. Using the Third Sylow Theorem
there are at least $p+1$ different Sylow $p$-subgroups in $G \leq
\GL_2(p)$, but there are exactly $p+1$ Sylow $p$-subgroups in
$\GL_2(p)$, so $G$ contains all the Sylow $p$-subgroups in
$\GL_2(p)$ and the subgroup they generate, thus $\SL_2(p) \leq G$.
\end{proof}

\begin{proposition}\label{cor:pknpkn_no_es}
For  a prime $p>3$ and an integer $n>1$, $\Z/p^n \times \Z/p^n$ is
not $\Ff$-Alperin in any
saturated fusion system $\Ff$ over a $p$-group $S$ where
$\Z/p^n\times \Z/p^n\lneq S$.
\end{proposition}
\begin{proof}
Consider $S$ a $p$-group and $\Z/p^n \times \Z/p^n \lneq S$. If we
assume that $\Z/p^n \times \Z/p^n$ is $\Ff$-radical then
$G\definicio \Aut_\Ff(\Z/p^n \times \Z/p^n)$ must be $p$-reduced,
and if $\Z/p^n \times \Z/p^n$ is $\Ff$-centric it is
self-centralizing in $S$, so taking the conjugation by an element
in $S \setminus \Z/p^n \times \Z/p^n$ we get that there exist an
element of order $p$ in $G$.

We can consider $\Aut(\Z/p^n \times \Z/p^n)$ as $2 \times 2$
matrices with coefficients in $\Z/p^n$ and with determinant non
divisible by $p$. In that case the reduction modulo $p$ induces a
short exact sequence: $$ \{1\} \rightarrow P \rightarrow
\Aut(\Z/p^n\times \Z/p^n) \stackrel{\rho}{\rightarrow} \GL_2(p)
\rightarrow \{1\} \, , $$ with $P$ a $p$-group. If the
intersection $P \cap G$ is not trivial, then we have a non trivial
normal $p$-subgroup in $G$ and $G$ would not be $p$-reduced, so
$G\cap P = \{ 1 \}$ and $\rho$ restricts to an injection of $G$ in
$\GL_2(p)$.

Using that $\rho(G)$ has an element of order $p$, and it is a
$p$-reduced subgroup of $\GL_2(p)$ and applying the Lemma
\ref{util} we get that $\rho(G)$ contains $\SL_2(p)$. In
particular we get that the matrix $\overline{A}\definicio
\big(\begin{smallmatrix} 1 & 1 \\ 0 & 1 \end{smallmatrix}\big)$ is
in $\rho(G)$, so we have an element in $A \in G \subset
\Aut(\Z/p^n \times \Z/p^n)$ which reduction modulo $p$ is $\overline
A$. $A$ must be matrix of the form $A=\big(\begin{smallmatrix}
1+ \xi p & 1+ \eta p \\ \lambda p & 1 + \mu p
\end{smallmatrix}\big)$, and an easy computation tell us that
$$
A^m \!\equiv\! \left(
\begin{array}{cc}
1\!+\!m \xi p\!+\!{m \choose 2} \lambda p \, & \,\,
m\!+\!{m \choose 2} \xi p \!+\! m \eta p \!+\! {m \choose 3} \lambda p \!+\! {m \choose 2} \mu p
\\m \lambda p & 1 \!+\! {m \choose 2} \lambda p \!+\! m \mu p
\end{array}
\right) \mod \! p^2.
$$
So, for $p>3$, we get $A^p\equiv \big(\begin{smallmatrix} 1 & p \\
0 & 1  \end{smallmatrix}\big) \, \mod p^2$  and as $n>1$, the
order of $A$ is bigger than $p$, which contradicts the fact that
$\rho$ is injective.
\end{proof}

If $p$ equals $3$ or $n$ equals $1$ then the thesis of the previous lemma is false, that is,
$\Z/p^n \times \Z/p^n$ could be $\Ff$-Alperin when $p=3$ or $n=1$. But then we can sharp our result in the following way:

\begin{lemma}\label{contiene31+2+}
If $p=3$ or $n=1$ and if
$\Z/p^n \times \Z/p^n$ is $\Ff$-Alperin
in a $p$-local finite group $(S,\Ff,\Ll)$  with
$\Z/p^n\times \Z/p^n\lneq S$, then $p^{1+2}_+ \leq S$.
\end{lemma}
\begin{proof}
Let $\Ff$ be a saturated fusion system over $S$, and suppose that
$P\definicio \Z/p^n \times \Z/p^n \lneq S$.
As in the proof of
Proposition \ref{cor:pknpkn_no_es}, if $P$ is $\Ff$-radical and $\Ff$-centric
then $\Aut_\Ff(P)$ contains $\SL_2(p)$.

Let $\pi\colon\Ll\longrightarrow\Ff^c$ be the centric
linking system, and take the
short exact sequence of groups induced by $\pi$:
$$
 {1} \rightarrow P \rightarrow \Aut_\Ll(P) \stackrel{\pi}{\rightarrow}
 \Aut_\Ff(P) \rightarrow {1}.
$$
Because $\Aut_\Ff(P)$ contains the special matrices over $\F_p$ we
have another short exact sequence:
$$
 {1} \rightarrow P \rightarrow M \stackrel{\pi}{\rightarrow} \SL_2(p)
\rightarrow {1}.
$$
To see that there are no more extensions that the split one,
consider the central subgroup~$V$ of $\SL_2(p)$ generated by the
involution $\big(\begin{smallmatrix} -1 & 0 \\ 0 &
-1\end{smallmatrix}\big)$. Since $p\geq 3$, multiplication by
$|V|=2$ is invertible in $P$ and so $H^k(V,P)=0$ for all $k>0$,
and because $V$ acts on $P$ without fixed points also
$H^0(V,P)=0$. Now the Hochschild-Serre spectral sequence
corresponding to the normal subgroup $V\leq \SL_2(p)$ give us that
$H^k(\SL_2(p),P)=0$ for all $k\geq 0$.

Now, as $H^2(\SL_2(p),P)=0$, the middle term $M$ of the short
exact sequence above equals $P\colon
\SL_2(p)\leq \Aut_\Ll(P)$.

As we are assuming in addition that the abelian $p$-group $P$ is
fully normalized then, as $\Ff$ is saturated, we obtain from
Definitions \ref{def:N} and \ref{def:centric_linking_system}(A)
that the Sylow $p$-subgroup of $\Aut_\Ll(P)$ is $N_S(P)$. So then
we have that $p^{1+2}_+\leq P:\Z/p \leq  N_S(P) \leq S$.
\end{proof}

Lemma \ref{util} also applies for giving some restrictions to
the family $C(p,r)$:

\begin{lemma}\label{lem:crfalperin}
Let $\Ff$ be a saturated fusion system over a $p$-group $S$ with
$p\geq 3$. If $C(p,r) \lneq S$ (with $r\geq 3$) is $\Ff$-centric and
$\Ff$-radical, then $\SL_2(p)\leq\Out_{\Ff}(C(p,r))\leq \GL_2(p)$.
\end{lemma}
\begin{proof}
If $C(p,r)$ is $\Ff$-radical, then $\Out_\Ff(C(p,r))$ must be
$p$-reduced. If we consider the projection $\rho$ from Lemma
\ref{lem:dpC(r)} we have that it must restrict to a monomorphism
in $\Out_\Ff(C(p,r))$ (otherwise we would have a non-trivial normal
$p$-subgroup) so we can consider \linebreak
$\Out_\Ff(C(p,r)) \leq \GL_2(p)$.

Now as $C(p,r)$ is $\Ff$-centric and different from $S$, we have an
element of order $p$ in $\Out_\Ff(C(p,r))$. Now, using again that
$\Out_\Ff(C(p,r))$ must be $p$ reduced, the result follows from
Lemma \ref{util}.
\end{proof}

%%%%%%%%%%%%%%%%%%%%%%%%%%%%%%%%%%%%%%%%%%%%%%%%%%%%%%%%%%%%%%%%
% 4. Non maximal class rang two $p$-groups.
%%%%%%%%%%%%%%%%%%%%%%%%%%%%%%%%%%%%%%%%%%%%%%%%%%%%%%%%%%%%%%%%

%\input{section4}
\section{Non-maximal class rank two $p$-groups}\label{metacyclic} \label{c(r)}

In this section we consider the non-maximal class rank two
$p$-groups for odd $p$, which are listed in the classification in
Theorem \ref{thm:clasificacionrango2}.

We begin with metacyclic groups:

\begin{theorem}\label{M-resistant}
Metacyclic $p$-groups are resistant for $p\geq 3$.
\end{theorem}
\begin{proof}
We prove that if $S$ is a metacyclic group then it cannot
contain any proper $\Ff$-Alperin subgroup.

Let $P$ be a proper subgroup of $S$, then it must be again
metacyclic.

If $P$ is not abelian, we can use Corollary \ref{cor:metanonabelian} to deduce that it cannot be
$\Ff$-Alperin.

If $P$ is abelian, using Corollaries \ref{cor:3.7} and
\ref{cor:3.6} it cannot be $\Ff$-Alperin unless $P\cong \Z/p^n
\times \Z/p^n$, but in this case $p^{1+2}_+$ should be a subgroup
of $M(p,r)$ by Lemma \ref{contiene31+2+}, which is impossible
because $p^{1+2}_+$ is not metacyclic.
\end{proof}

In the study of $C(p,r)$ in this section we assume $r\geq4$: for $r=3$ we have that $C(p,3)\!\cong\!
p^{1+2}_+$, which it is a maximal nilpotency class $p$-group and the fusion systems over that
group are studied in \cite{RV}.

\begin{theorem}\label{C-resistant}
If $r>3$ and $p\geq 3$ then $C(p,r)$ is resistant.
\end{theorem}
\begin{proof}
Let $\Ff$ be a saturated fusion system over $C(p,r)$. The possible
proper $\Ff$-Alperin subgroups are proper $p$-centric subgroups,
and using Lemma \ref{lem:pcentricr} these are isomorphic to
$\Z/p^{r-2}\times\Z/p$. But as $r>3$, $\Z/p^{r-2}\times\Z/p$
cannot be $\Ff$-radical by Corollary \ref{cor:3.6}.
\end{proof}

It remains to study the non-maximal nilpotency class groups of
type $G(p,r;\epsilon)$. Remark~\ref{rmk:idBiG} tells us that in
this section we have to consider all of them but $G(3,4;\pm1)$.

The study of groups of type $G(p,r;\epsilon)$ is divided in two cases: for $p\geq5$ these are resistant groups, while for $p=3$ we
obtain saturated fusion systems with proper $\Ff$-Alperin
subgroups.

\begin{theorem}\label{Thm:g(r,e)isresistant}
If $p> 3$ and $r\geq 4$, $G(p,r;\epsilon)$ is resistant.
\end{theorem}

The proof needs the following lemma:

\begin{lemma}\label{aut(c(r-1))}
Let $\Ff$ be a saturated fusion system over $G(p,r;\epsilon)$ with $p>3$
and $r\geq4$. Then $C(p,r-1)<G(p,r;\epsilon)$ is not $\Ff$-radical. Moreover,
$\Aut_{\Ff}(C(p,r-1))$ is a subgroup of the lower triangular
matrices with first diagonal entry $\pm 1$.
\end{lemma}
\begin{proof}
Consider $p\geq5$ and
assume that $C(p,r-1)$ is $\Ff$-radical, then by Lemma \ref{lem:crfalperin}
$\SL_2(p)\leq\Out_{\Ff}(C(p,r-1))$ and therefore the matrix
$\big(\begin{smallmatrix} x & 0 \\ 0 & x^{-1}
\end{smallmatrix}\big)$, where $x$ is a primitive $(p-1)$-th root
of the unity in $\F_p$, is the image of some $\varphi\in
\Aut_{\Ff}(C(p,r-1))$ by the composition
$$
\xymatrix{
\Aut_{\Ff}(C(p,r-1)) \ar[r]^{\pi} & \Out_{\Ff}(C(p,r-1)) \ar[r]^{\qquad\rho} & \GL_2(p)
}
$$
of the projection and $\rho$ from Lemma
\ref{lem:dpC(r)}. By Definition \ref{def:N} this morphism $\varphi$ must lift
to $\Aut_{\Ff}(G(p,r;\epsilon))$ as a morphism that maps $a$ to $a^x$ (up
to $b^nc^{mp}$ multiplication). But automorphisms of $G(p,r;\epsilon)$ map
$a$ to $a^{\pm 1}$ (up to $b^nc^{mp}$ multiplication), and $-1$ is
not a primitive $(p-1)$-th root of the unity in $\F_p$ for $p\geq5$.
Therefore $C(p,r-1)$ is not $\Ff$-radical.

Finally, as $C(p,r-1)$ is not $\Ff$-radical, then
$\pi(\Aut_{\Ff}(C(p,r-1)))<N_{\GL_2(p)}(V)$ where $V$ is the group
generated by $\big(\begin{smallmatrix} 1 & 0 \\ -1 & 1
\end{smallmatrix}\big)$, that is, the projection $\pi(\Aut_{\Ff}(C(p,r-1)))$ is a
subgroup of the lower triangular matrices. Then all morphisms in
$\Aut_{\Ff}(C(p,r-1))$ must lift to $\Aut_{\Ff}(G(p,r;\epsilon))$, and again
only those with $\pm 1$ in the first diagonal entry are allowed,
what proves the second part of the lemma.
\end{proof}

\begin{proof}[Proof of Theorem \ref{Thm:g(r,e)isresistant}]
If $r>4$ we have just to consider the case of $C(p,r-1)$ being
$\Ff$-radical, and Lemma \ref{aut(c(r-1))} shows that it cannot
be.

If $r=4$ it remains to check what happens with the rank two
elementary abelian subgroups in $G(p,4;\epsilon)$. According to
Lemma \ref{pcentricenG(re)} there are exactly $p$ of these
subgroups, namely, $V_i \definicio  \langle ab^i,c^p \rangle $ for
$i=0,...,p-1$, and all of them lie inside $C(p,3)\cong p^{1+2}_+$.
Notice that conjugation by $c$ permutes all of them cyclically and
so they are $\Ff$-conjugate in any saturated fusion system over
$G(p,4;\varepsilon)$. Thus if any of them is $\Ff$-radical then
all of them are $\Ff$-radical.

If this is the case, fix $x$ a primitive $(p-1)$-th root of the
unity in $\F_p$ and let $V_i$ be one of these rank two
$p$-subgroups with $\F_p$ basis $\langle ab^i,c^{-ep} \rangle$.
Then $\SL_2(p)\leq\Aut_{\Ff}(V_i)$ by Lemma \ref{util}. Now, the
element $\big(\begin{smallmatrix} x & 0
\\ 0 & x^{-1}\end{smallmatrix}\big)\in\Aut_{\Ff}(V)$ must lift to
an automorphism of $p^{1+2}_+$ by Definition~\ref{def:N}. The image of this
extension in $\Out_\Ff(p^{1+2}_+)=\GL_2(p)$ is a matrix $L_i$ with
$x$ as eigenvalue and with determinant $x^{-1}$, so it has
$x^{-2}$ as the other eigenvalue. Notice that each $V_i$ gives
a different matrix $L_i \in \GL_2(p)$, so different elements  of
order $p-1$ in $\Out_{\Ff}(p^{1+2}_+)$. But the description of
$\Aut_{\Ff}(p^{1+2}_+)$ in Lemma \ref{aut(c(r-1))} shows us that
there are not such matrices $L_i$ if $p>5$ and at most one when $p=5$
(cf. \cite[Section 4]{RV}).
\end{proof}

So it remains to check the cases with $p=3$, and as this section
only copes with the non maximal nilpotency class groups, we consider
$r \geq 5$.

\begin{lemma}\label{lem:accionesenCpr}
Fix $r\geq4$ and let $W$ be either $\SL_2(3)$ or $\GL_2(3)$.
Given a faithful representation of $\Z/3$ in $\Out(C(p,r-1))$ and
any group $S$ fitting in the exact sequence
$$
\xymatrix{
1 \ar[r]         &  C(3,r-1) \ar[r]  & S \ar[r]  &   \Z/3 \ar[r] & 1 \\
}
$$
that induces the given representation of $\Z/3$,
there is a finite group $G$ which fits in the following commutative diagram
$$
\xymatrix{
1 \ar[r]         &  C(3,r-1) \ar[r] \ar@{=}[d] & S \ar[r] \ar[d] &   \Z/3 \ar[r] \ar[d] & 1 \\
1 \ar[r]         &  C(3,r-1) \ar[r]        & G \ar[r]        & W     \ar[r]   & 1 \\
}
$$
where the last column is an inclusion of the Sylow $3$-subgroup in $W$.
\end{lemma}
\begin{proof}
% If we consider a trivial representation of $\Z/3$ in $\Out(C(p,r-1))$, we can take
% $G$ as the direct product $C(p,r-1)\times W$.

Fix a faithful representation of $\Z/3$ in $\Out(C(p,r-1))$.

Consider first $r=4$ or $W=\GL_2(3)$.
Using the study of the group $\Aut(C(3,r-1))$ in Lemma \ref{lem:dpC(r)}
we have that there is only one faithful representation
of $W$ in $\Out(C(3,r-1))$ up to conjugation.

To compute the possible equivalence class of extensions in the exact sequence
$1 \to C(3,r-1) \to G \to W \to 1$ with the fixed representation of $W$ in $\Out(C(3,r-1))$
we have to compute $H^2(W;Z(C(3,r-1)))$
\cite[Theorem IV.6.6]{brown}, getting $H^2(W;\Z/3^{r-3})\cong \Z/3$.

Consider now $\Z/3$, the Sylow $3$-subgroup of $W$, and consider
the following diagram:
$$
\xymatrix{
1 \ar[r]         &  C(3,r-1) \ar[r] \ar@{=}[d] & S \ar[r] \ar[d] &   \Z/3 \ar[r] \ar[d] & 1 \\
1 \ar[r]         &  C(3,r-1) \ar[r]        & G \ar[r]        & W     \ar[r]   & 1 \\
}
$$
with exact rows and do the same computation for the first row
in cohomology. This gives us that there are three possible equivalence
classes of extensions $S$
fitting in that exact sequence. Moreover the map induced in
cohomology $H^2(W;\Z/3^{r-3}) \to H^2(\Z/3;\Z/3^{r-3})$ is
an isomorphism (use a transfer argument), so we have that all
three $S$ appear as the Sylow $3$-subgroup of $G$.

Consider now the case $r\geq5$ and $W=\SL_2(3)$. Now by Lemma \ref{lem:dpC(r)}
we have two possible faithful representations of $W$ in $\Out(C(p,r-1))$ up
to conjugation. Now we do the same computations and we use the same argument
as before for each action, getting $6$ possible equivalence classes of
extensions which appear as
the Sylow $3$-subgroup of the $6$ possible $G$'s.
\end{proof}

\begin{remark}
The previous assertion fails for $p>3$ because in the cohomology groups computation
we get that $H^2(\SL_2(p);Z(C(p,r-1)))$ and
$H^2(\GL_2(p);Z(C(p,r-1)))$ are trivial, so the only extensions are
the split ones. So in both cases the Sylow $p$-subgroup has $p$-rank three.
\end{remark}

\begin{remark}\label{rmk:accionesenCpr}
We are interested in the possible extensions $1\to C(3,r-1)\to S \to \Z/3 \to 1$
such that the $3$-rank of $S$ is two. Looking at the
groups appearing in the classification
in Theorem \ref{thm:clasificacionrango2}
we get that the only possible ones are $C(3,r)$, $G(3,r;1)$ and $G(3,r;-1)$.
\end{remark}

Fix the notation $C(3,r-1)\cdot_\epsilon W$ for the groups in Lemma
\ref{lem:accionesenCpr} with Sylow $3$-subgroup isomorphic to $G(3,r;\epsilon)$.

\begin{theorem}\label{sfs-G-3}
The $3$-local finite groups over $G(3,r;\epsilon)$ with
$r \geq 5$ and at least one proper $\Ff$-Alperin subgroup
are classified by the following parameters:
\begin{center}
\begin{table}[H]
$$
\begin{array}{|c|c|c|}
\hline \pb \Out_\Ff(G(3,r;\epsilon)) &  \Out_\Ff(C(3,r-1)) &\mbox{Group} \\
\hline \hline
\pb \Z/2 & \SL_2(3) & C(3,r-1)\cdot_\epsilon\SL_2(3) \\
\hline
\pb \Z/2\times\Z/2 & \GL_2(3) & C(3,r-1)\cdot_\epsilon\GL_2(3) \\
\hline
\end{array}
$$
\caption{s.f.s. over $G(3,r;\epsilon)$ for $r\geq5$.}
\end{table}
\end{center}
where the second column gives the outer automorphism group over
the only proper $\Ff$-Alperin subgroup, which is isomorphic to
$C(3,r-1)$. In the last column we refer to non split extensions
such that the Sylow $3$-subgroup is isomorphic to
$G(3,r;\epsilon)$ for $\epsilon = \pm1$ and it induces the desired
fusion system.
\end{theorem}
\begin{proof}
Assume now that $C(3,r-1) \leq G(3,r;\epsilon)$ is $\Ff$-Alperin. Then
$\Out_\Ff(C(3,r-1)) \cong \SL_2(3)$ or
$\Out_\Ff(C(3,r-1)) \cong \GL_2(3)$.

If we look at the possible $\Out_\Ff(G(3,r;\epsilon))$ we get that
it is $\Z/2$, generated by a matrix which induces $-\Id$ in the
identification $\Out_\Ff(C(3,r-1)) \leq \GL_2(3)$ or $\Z/2 \times
\Z/2$, generated by $-\Id$ and a matrix of determinant $-1$ in
$\Out_\Ff(C(3,r-1)) \leq \GL_2(3)$. From here we deduce that if
$\Out_\Ff(G(3,r;\epsilon))\cong\Z/2$ and $C(3,r-1)$ is
$\Ff$-radical then $\Out_\Ff(C(3,r-1))\cong \SL_2(3)$,
while if $\Out_{\Ff}(G(3,r;\epsilon))\cong\Z/2\times \Z/2$ then
$\Out_{\Ff}(C(3,r-1))\cong\GL_2(3)$ what completes the table.

All of them are saturated because are the fusion systems of the
finite groups described in Lemma \ref{lem:accionesenCpr}.
\end{proof}

%%%%%%%%%%%%%%%%%%%%%%%%%%%%%%%%%%%%%%%%%%%%%%%%%%%%%%%%%%%%%%%%
% 5. Maximal class rang two $p$-groups.
%%%%%%%%%%%%%%%%%%%%%%%%%%%%%%%%%%%%%%%%%%%%%%%%%%%%%%%%%%%%%%%%

%\input{section5}
\section{Maximal nilpotency class rank two $p$-groups}\label{maximal_class}

In this section we classify the $p$-local finite groups over
$p$-groups of maximal nilpotency class and $p$-rank two. Recall
that by Corollary \ref{cor:rank2} we have just to classify the
saturated fusion systems over these groups.

Consider $S$ a $p$-rank two maximal nilpotency class $p$-group of order $p^r$.
For $r=2$ then $S\cong\Z/p\times\Z/p$, which is resistant by
Corollary \ref{cor:abelian}.
If $r=3$ then $S\cong p^{1+2}_+$ and this case has been studied in
\cite{RV}.
For $r\geq4$ all the $p$-rank two maximal nilpotency class groups
appear only at $p=3$, and we
use the description and properties given in Appendix \ref{rank-two-p-groups}.

The description of the maximal nilpotency class $3$-groups of
order bigger than $3^3$ depends on three parameters $\beta$,
$\gamma$ and $\delta$, and we use the notation given in Theorem
\ref{thm:maxclass}, so we call the groups
$B(3,r;\beta,\gamma,\delta)$ and $\{s,s_1,s_2, \dots, s_{r-1}\}$
the generators.

First we consider the non split case, that is, $\delta>0$.

\begin{theorem}\label{B(n;0,b,c,d)resistant}
Every group of type $B(3,r;\beta,\gamma,\delta)$,
$\delta>0$, is resistant.
\end{theorem}
\begin{proof}
Let $\Ff$ be a saturated fusion system over
$B(3,r;\beta,\gamma,\delta)$. Using Alperin's fusion theorem for saturated fusion systems
 (Theorem \ref{teo:alperin}) it is enough to see whether
$B(3,r;\beta,\gamma,\delta)$ is the only $\Ff$-Alperin subgroup.
So let $P$ be a proper subgroup.

If $P$ has $3$-rank one then, it is cyclic and we can apply
Corollary \ref{cor:3.7} to obtain that $P$ cannot be
$\Ff$-Alperin.

Now assume that $P$ has $3$-rank
two. Then, by Theorem \ref{thm:clasificacionrango2}, $P$ is one of the following:
\begin{itemize}
\item $M(3,r)$ non abelian: according to Corollary \ref{cor:metanonabelian} $P$ can not be $\Ff$-Alperin.
\item $G(3,r;\epsilon)$ group: it cannot be $\Ff$-Alperin by
Corollary \ref{cor:Grenopuede}.
\item $B(3,m;\beta,\gamma,\delta)$ with $m<r$: by Lemma \ref{B(n)_noes} $P$
cannot be $\Ff$-Alperin.
\item $C(3,r)$ group: as $3^{1+2}_+=C(3,3)$ is contained in $C(3,r)$ we obtain that
$3^{1+2}_+\leq B(3,r;\beta,\gamma,\delta)$.
\item Abelian: say $P=\Z/3^m \times \Z/3^n$. If $m\neq n$ then $P$ cannot be
$\Ff$-Alperin by Corollary~\ref{cor:3.6}, and if $m=n$ then again
$3^{1+2}_+\leq B(3,r;\beta,\gamma,\delta)$ by Lemma \ref{contiene31+2+}.
\end{itemize}
So if $P$ is $\Ff$-Alperin then $B(3,r;\beta,\gamma,\delta)$ must
contain $3^{1+2}_+$.

We finish the proof by showing that $3^{1+2}_+\nleq B(3,r;\beta,\gamma,\delta)$
since we are in the non split case ($\delta>0$).
Consider the short exact
sequence:
$$
 {1} \rightarrow \gamma_1 \rightarrow B(3,r;\beta,\gamma,\delta) \stackrel{\pi}\rightarrow
 \Z/3 \rightarrow {1}.
$$ If $\pi(3^{1+2}_+)$ is trivial then $3^{1+2}_+\leq \gamma_1$.
But by \cite[Satz III.{\S}14.17]{huppert} $\gamma_1$ is metacyclic, and
consequently all its subgroups are metacyclic too. Thus $\gamma_1$
cannot contain the $C(3,3)$ group. We obtain then the short exact
sequence: $$
 {1} \rightarrow \Z/3 \times \Z/3 \rightarrow 3^{1+2}_+ \stackrel{\pi}\rightarrow
 \Z/3 \rightarrow {1}.
$$
As this sequence splits, the same would holds for the exact
sequence involving  \\
$B(3,r;\beta,\gamma,\delta)$, and this is not the case.
\end{proof}

We now consider the split case, that is $\delta=0$. We prove
that for $\beta=1$ these maximal nilpotency class groups are resistant,
while for $\beta=0$ they are not.
We use the information about $B(3,r;\beta,\gamma,0)$ contained in the appendix.
According to Alperin's fusion theorem for saturated fusion systems (Theorem \ref{teo:alperin})
we must focus on the $\Ff$-Alperin subgroups. The next two lemmas list the
subgroups candidates for being $\Ff$-Alperin in a saturated fusion system
over $B(3,r;0,\gamma,0)$ and $B(3,r;1,0,0)$:

\begin{lemma}\label{B(n;0,0,c,0)Alperin}
Let $\Ff$ be a saturated fusion system over $B(3,r;0,\gamma,0)$
and let $P$ be a proper $\Ff$-Alperin subgroup. Then $P$ is one of
the following table:
\begin{center}
\begin{tabular}{|@{}c@{}|l|}
\hline
\begin{tabular}{p{23mm}|p{29mm}}
 Isomorphism type & Subgroup (up to conjugation)\\
\end{tabular} & Conditions \\
\hline \hline
\begin{tabular}{p{23mm}|p{29mm}}
\pb $\Z/3^k \times \Z/3^k$ & $\gamma_1=\langle s_1,s_2 \rangle$
\end{tabular}
 & $r=2k+1$.\\
\hline
\begin{tabular}{p{23mm}|p{29mm}}
\pbbb $3^{1+2}_+$ & $E_i \!\definicio\! \langle \zeta,\zeta',s {s_1}^i \rangle$
\\ \hline
\pbbb $\Z/3 \times \Z/3$ & $V_i \definicio \langle \zeta,s {s_1}^i \rangle$
\end{tabular}
&  \parbox{61mm}{$\zeta\!=\!{s_2}^{3^{k-1}}\!, \zeta'\!=\!{s_1}^{3^{k-1}}$ for $r\!=\!2k+1$, \\
   $\zeta\!=\!{s_1}^{3^{k-1}},  \zeta'\!=\!{s_2}^{-3^{k-2}}$ for $r\!=\!2k$, \\
   $i\in\{-1,0,1\}$ if $\gamma=0$ and  \\
   $i=0$ if $\gamma=1,2$.}
\\
\hline
\end{tabular}
%\end{tabular}
\end{center}
Moreover, for the subgroups in the table, $\gamma_1$ and any $E_i$
are always $\Ff$-centric, while any $V_i$  is $\Ff$-centric only
if it is not $\Ff$-conjugate to $\langle
\zeta,\zeta'\rangle\cong\Z/3\times\Z/3$.
\end{lemma}
\begin{proof}
If $P$ is $\Ff$-Alperin then arguing as in Theorem \ref{B(n;0,b,c,d)resistant} we obtain
that $P\cong\Z/3^n \times \Z/3^n$ or $P\cong C(3,n)$. Then, using
Lemma \ref{B(n;0,0,c,0)centricos} we reach the subgroups in the statement.

To check that indeed $\gamma_1$ is $\Ff$-centric it is enough to
notice that it is self-centralizing and that $\gamma_1$ is
$\Ff$-conjugate just to itself in any fusion system $\Ff$
($\gamma_1$ is strongly characteristic in $B(3,r;0,\gamma,0)$). It
is clear that the copies of $3^{1+2}_+$ are $\Ff$-centric because
these are the only copies lying in $B(3,r;0,\gamma,0)$, and the
center of all of them is $\langle \zeta \rangle$. To conclude the
lemma, notice that the only copies of $\Z/3\times\Z/3$ in
$B(3,r;0,\gamma,0)$ are the $V_i$'s and $\langle \zeta,\zeta'
\rangle$, and that the former are self-centralizing while the
centralizer of the latter is $\gamma_1$.
\end{proof}

\begin{lemma}\label{B(n;0,1,0,0)Alperin}
Let $\Ff$ be a saturated fusion system over $B(3,r;1,0,0)$
and $P$ be a proper $\Ff$-Alperin subgroup. Then $P$ is one of
the following table:
\begin{center}
\begin{tabular}{|@{}c@{}|l|}
\hline
\begin{tabular}{p{29mm}|p{26mm}}
Isomorphism type & Subgroup (up to conjugation) \\
\end{tabular} & Conditions\\
\hline \hline
\begin{tabular}{p{29mm}|p{26mm}}
\pbb $\Z/3^{k-1} \times \Z/3^{k-1}$ & $\gamma_2=\langle s_2,s_3 \rangle$
\end{tabular} & $r=2k$.\\
\hline
\begin{tabular}{p{29mm}|p{26mm}}
\pbbb $3^{1+2}_+$ & $E_0\definicio \langle \zeta,\zeta',s \rangle$ \\
\hline
\pbb $\Z/3 \times \Z/3$ & $V_0 \definicio \langle \zeta,s \rangle$ \\
\end{tabular} & \parbox{68mm}{$\zeta={s_2}^{3^{k-1}}$, $\zeta'={s_3}^{-3^{k-2}}$
 for $r=2k+1$, \\
   $\zeta={s_3}^{3^{k-2}}$, $\zeta'={s_2}^{3^{k-2}}$ for $r=2k$.}\\
\hline
\end{tabular}
\end{center}
Moreover, for the subgroups in the table, $\gamma_2$ and $E_0$
are always $\Ff$-centric, and $V_0$ is $\Ff$-centric only if it is
not $\Ff$-conjugate to $\langle
\zeta,\zeta'\rangle\cong\Z/3\times\Z/3$.
\end{lemma}
\begin{proof}
The reasoning is totally analogous to that of the previous lemma using Lemma \ref{B(n;0,1,0,0)centricos}.
\end{proof}

\begin{remark}\label{rmk:bases}
The isomorphism between $3^{1+2}_+$ and $E_i$ is given by
$a\mapsto \zeta'$, $b\mapsto s{s_1}^i$ and $c\mapsto \zeta$, where
$a$, $b$ and $c$ are the generators of $3^{1+2}_+=C(3,3)$ given in
Theorem \ref{thm:clasificacionrango2}. Therefore the morphisms in
$\Out(E_i)$ are described as in \cite[Lemma 3.1]{RV} by means of
the mentioned isomorphism, and for $\Out(V_i)$ we choose the
ordered basis $\{\zeta,s{s_1}^i\}$.

When studying a saturated fusion system $\Ff$ over $B(3,r;0,\gamma,0)$ or $B(3,r;1,0,0)$ it is
enough to study the representatives subgroups given in the tables of Lemmas
\ref{B(n;0,0,c,0)centricos} and \ref{B(n;0,1,0,0)centricos} because $\Ff$-properties are invariant
under conjugation.
\end{remark}

As a last step before the classification itself, we work out some
information on lifts and restrictions of automorphism of some subgroups
of $B(3,r;\beta,\gamma,0)$. Given a saturated fusion system $\Ff$
over $B(3,r;\beta,\gamma,0)$ subsequently it will be advantageous
to consider for every subgroup of rank two $P\leq
B(3,r;\beta,\gamma,0)$ the Frattini map $
\Out(P)\stackrel{\rho}\rightarrow \GL_2(3)$, which kernel is a $3$
group, and its restriction $$
    \Out_{\Ff}(P)\stackrel{\rho}\rightarrow \GL_2(3) \,.
$$ Notice that by Remark \ref{rmk:orderW} this restriction is a
monomorphism for $P=B(3,r;\beta,\gamma,0)$. For $P=\gamma_1$ or
$\gamma_2$ it is a monomorphism if $P$ is $\Ff$-Alperin as
$\Out_{\Ff}(P)$ is $3$-reduced. For $P=E_i$, $V_i$ they are
inclusions by \cite[Lemma 3.1]{RV} and by definition respectively.
In the case this restriction is a monomorphism we identify
$\Out_{\Ff}(P)$ with its image in $\GL_2(3)$ without explicit
mention. We divide the results in three lemmas, the first gives a
description of $\Out_{\Ff}(P)$ for~$P$ an $\Ff$-Alperin subgroup,
the second copes with restrictions and the third with lifts.
\begin{lemma}\label{B(n)F-automorphisms}
Let $\Ff$ be a saturated fusion system over $B(3,r;\beta,\gamma,0)$.
Then:
\begin{enumerate}[\rm (a)]
\item $\rho(\Out_{\Ff}(B(3,r;\beta,\gamma,0)))$ is a subgroup of the lower
triangular matrices,
\item $\Out_{\Ff}(B(3,2k;0,\gamma,0))\leq\Z/2\times\Z/2$,
\item $\Out_{\Ff}(B(3,2k+1;0,1,0))\leq\Z/2$,
\item $\Out_{\Ff}(B(3,r;1,0,0))\leq\Z/2$ and
\item if $P=\gamma_1$, $\gamma_2$, $E_i$ or $V_i$ is
$\Ff$-Alperin then $\Out_{\Ff}(P)=\SL_2(3)$ or $\GL_2(3)$.
\end{enumerate}
\end{lemma}
\begin{proof}
For the first claim notice that from the order of
$\Aut(B(3,r;\beta,\gamma,0))$ we deduce that a $3'$ element $\varphi$ should
have order $2$ or $4$. Now, as the Frattini subgroup of $B(3,r;\beta,\gamma,0)$ is $\langle s_2,s_3\rangle$,
the projection on the Frattini quotient becomes
$$
    \Out_{\Ff}(B(3,r;\beta,\gamma,0))\stackrel{\rho}\rightarrow \GL_2(3)
$$
where if $\overline{\varphi}$ maps $s$ to $s^e {s_1}^{e'} {s_2}^{e''}$ and $s_1$ to $s_1^{f'}
s_2^{f''}$ then $\rho(\overline{\varphi})=\big(\begin{smallmatrix} e & 0\\ e' &
f'\end{smallmatrix}\big)$. So we obtain lower triangular matrices. Checking cases leads to obtain
that the order of $\varphi$ is two, and that $\varphi \in \{ \big(\begin{smallmatrix} 1 & 0\\ 0 &
-1\end{smallmatrix}\big), \big(\begin{smallmatrix} -1 & 0\\ 0 & 1\end{smallmatrix}\big),
\big(\begin{smallmatrix} -1 & 0\\ 0 & -1\end{smallmatrix}\big), \big(\begin{smallmatrix} 1 & 0\\ 1
& -1\end{smallmatrix}\big), \big(\begin{smallmatrix} -1 & 0\\ 1 & 1\end{smallmatrix}\big),$
$\big(\begin{smallmatrix} 1 & 0\\ -1 & -1\end{smallmatrix}\big), \big(\begin{smallmatrix} -1 & 0\\
-1 & 1\end{smallmatrix}\big) \}$. For $\beta=0$, $\gamma=1$ and odd $r$ or $\beta=1$ Lemma
\ref{B(n)automorphisms} tells us that $e$ must be equal to $1$, and then it is easily deduced that
$\Out_{\Ff}(B(3,r;\beta,\gamma,0))$ must have order two.

For the last point just use Lemma \ref{util} and that
$[\GL_2(3):\SL_2(3)]=2$.
\end{proof}

Now we focus on the restrictions. For the study of the possible
saturated fusion systems $\Ff$  %and the structure described in \cite{BB},
it is enough to consider diagonal matrices of
$\Out_{\Ff}(B(3,r;\beta,\gamma,0))$ instead of lower triangular
ones. This is so because every $\Z/2$ and $\Z/2\times\Z/2$ in
$\Out_{\Ff}(B(3,r;\beta,\gamma,0))$ is
$B(3,r;\beta,\gamma,0)$-conjugate to a diagonal one.

\begin{lemma}\label{B(n)restrictions}
Let $\Ff$ be a saturated fusion system over
$B(3,r;\beta,\gamma,0)$. Then:
\begin{enumerate}[\rm (a)]
\item  If $\beta=0$ the restrictions  of
 the elements $\overline{\varphi}\in
\Out_{\Ff}(B(3,r;0,\gamma,0))$ to $\Out_{\Ff}(\gamma_1)$ are given by the following table,
where it is also described the permutation of the $E_i$'s induced
by~$\overline{\varphi}$, and the restrictions to
$\Out_{\Ff}(E_{i_0})$ for $i_0\in \{-1,0,1\}$ such that $E_{i_0}$
is fixed by the permutation.
%some $\psi\in\overline{\varphi}$.
$$
\begin{array}{|c|c|c|c|}
\hline
\pbb
 \Out_{\Ff}(B(3,r;0,\gamma,0))&\big(\begin{smallmatrix} 1 & 0\\
0 & -1\end{smallmatrix}\big)&\big(\begin{smallmatrix} -1 & 0\\ 0 &
1\end{smallmatrix}\big) &\big(\begin{smallmatrix} -1 & 0\\ 0 &
-1\end{smallmatrix}\big)\\
\hline %%
\pbb
 \Out_{\Ff}(\gamma_1) & \big(\begin{smallmatrix} -1 & 0\\ 0 &
-1\end{smallmatrix}\big)& \big(\begin{smallmatrix} 1 & 0\\ 0 &
-1\end{smallmatrix}\big)&\big(\begin{smallmatrix} -1 & 0\\ 0 &
1\end{smallmatrix}\big)\\
\hline %%
\pb $\Ff$\mbox{-conjugation} & E_1\leftrightarrow E_{-1} & E_1\leftrightarrow E_{-1} & -\\
\hline
\pbb
\Out_{\Ff}(E_{i_0})\text{, $r$ odd} &
\big(\begin{smallmatrix} -1 & 0\\ 0 & 1\end{smallmatrix}\big) &
\big(\begin{smallmatrix} 1 & 0\\ 0 & -1\end{smallmatrix}\big) &
\big(\begin{smallmatrix} -1 & 0\\ 0 & -1\end{smallmatrix}\big)\\
\hline
\pbb
\Out_{\Ff}(E_{i_0})\text{, $r$ even} &
\big(\begin{smallmatrix} -1 & 0\\ 0& 1\end{smallmatrix}\big) &
\big(\begin{smallmatrix} -1 & 0\\ 0& -1\end{smallmatrix}\big) &
\big(\begin{smallmatrix} 1 & 0\\ 0& -1\end{smallmatrix}\big) \\
\hline
\end{array}
$$

\item  If $\beta=1$ (thus $\gamma=0$)
the restrictions to $\Out_{\Ff}(\gamma_2)$ and to $\Out_{\Ff}(E_0)$ of outer automorphisms
$\overline{\varphi}\in \Out_{\Ff}(B(3,r;1,0,0))$ are given by the
following table.
$$
\begin{array}{|c|c|}
\hline \pbb
 \Out_{\Ff}(B(3,r;1,0,0))&\big(\begin{smallmatrix} 1 & 0\\
0 & -1\end{smallmatrix}\big)\\
\hline %%
\pbb
 \Out_{\Ff}(\gamma_2) & \big(\begin{smallmatrix} -1 & 0\\ 0 &
-1\end{smallmatrix}\big)\\
\hline %%
\pbb
\Out_{\Ff}(E_0) & \big(\begin{smallmatrix} -1 & 0\\ 0&
1\end{smallmatrix}\big) \\ \hline
\end{array}
$$

\item For the outer automorphism groups of $E_i$ and
$V_i$ we have the following restrictions:
$$
\begin{array}{|c|c|c|c|}
\hline \pbb \Out_{\Ff}(E_i) &
\big(\begin{smallmatrix} 1 & 0\\ 0 &-1\end{smallmatrix}\big) &
\big(\begin{smallmatrix} -1 & 0\\ 0 &1\end{smallmatrix}\big) &
\big(\begin{smallmatrix} -1 & 0\\ 0 & -1\end{smallmatrix}\big)\\
\hline \pbb \Out_{\Ff}(V_i) &
\big(\begin{smallmatrix}-1&0\\0&-1\end{smallmatrix}\big) &
\big(\begin{smallmatrix} -1 & 0\\ 0 &1\end{smallmatrix}\big) &
\big(\begin{smallmatrix} 1 & 0\\0 &-1\end{smallmatrix}\big) \\
\hline
\end{array}
$$
\end{enumerate}
\end{lemma}
\begin{proof}
These are the only possible elements by Lemma
\ref{B(n)F-automorphisms}, and the restrictions are computed
directly using the explicit form of the subgroups given in Lemmas
\ref{B(n;0,0,c,0)Alperin} and \ref{B(n;0,1,0,0)Alperin} using the
basis described in Remark \ref{rmk:bases}.
\end{proof}

Finally we reach the lemma about lifts:
\begin{lemma}\label{B(n)liftings}
Let $\Ff$ be a saturated fusion system over
$B(3,r;\beta,\gamma,0)$, and $P$ be one of the proper subgroups appearing in
the tables of Lemmas \ref{B(n;0,0,c,0)Alperin} or \ref{B(n;0,1,0,0)Alperin}.
Then every diagonal outer automorphism of $P$ in $\Ff$
(like those appearing in Lemma \ref{B(n)restrictions}) can be
lifted to the whole $B(3,r;\beta,\gamma,0)$. In particular, every admissible diagonal outer
automorphism of $P$ in $\Ff$ is listed in the tables of restrictions of Lemma
\ref{B(n)restrictions}.
\end{lemma}
\begin{proof}
We study the two cases $\beta=0$ and $\beta=1$ separately.

If $\beta=0$ we begin with the morphisms from $B(3,r;0,\gamma,0)$
restricted to $\gamma_1$. If $\gamma_1$ is not $\Ff$-Alperin then
every morphism in $\Aut_{\Ff}(\gamma_1)$ can be lifted by Theorem
\ref{teo:alperin}. Suppose then that $\gamma_1$ is $\Ff$-Alperin.
Take $\varphi\in \Out_{\Ff}(\gamma_1)$ appearing in the table of
Lemma \ref{B(n)restrictions} and consider the images $c_s,c_{s^2}\in
\Out_{\Ff}(\gamma_1)$ of the restrictions to $\gamma_1$ of
conjugation by $s$ and $s^2$. Then it can be checked that $\varphi
c_s {\varphi}^{-1}$ equals $c_s$ or $c_{s^2}$. Now apply (II) from
Definition \ref{def:N}.

For $E_i\leq B(3,r;0,\gamma,0)$ there is a little bit more work to
do. If $E_i$ is not $\Ff$-Alperin apply Theorem \ref{teo:alperin}
again. So suppose $E_i$ is $\Ff$-Alperin, take
$\overline{\varphi}\in\Out_\Ff(E_i)$ and compute $N_\varphi$ from
Definition \ref{def:N}. The normalizer of $E_i$ is $\langle
\zeta,\zeta'',sa \rangle\cong (\Z/9\times\Z/3):\Z/3$ where the
order of~$\zeta''$ equals $9$, and $\varphi c_{\zeta''}
{\varphi}^{-1}$ equals $c_{\zeta''}$ or $c_{{\zeta''}^2}$. So
$\varphi$ can be lifted to the normalizer (Definition~\ref{def:N})
and, as this subgroup cannot be $\Ff$-Alperin by Lemma
\ref{B(n;0,0,c,0)Alperin}, we can extend again to the whole
$B(3,r;0,\gamma,0)$ by Theorem \ref{teo:alperin}. The last case is
to consider $V_i\leq E_i$. The details are analogous with
$c_{\zeta'}$.

If $\beta=1$, then $P$ can be identified with a centric subgroup of
$B(3,r-1;0,0,0)<B(3,r;1,0,0)$ (see proof of Lemma
\ref{B(n;0,1,0,0)centricos}) and then use the arguments above to
extend the morphisms to $B(3,r-1;0,0,0)$. We then use saturateness
to extend the morphism to $B(3,r;1,0,0)$.
\end{proof}

\begin{theorem}\label{Thm:B(r;0,1,0,0)resistant}
Every rank two $3$-group isomorphic to $B(3,r;1,0,0)$ is resistant.
\end{theorem}
\begin{proof}
Using Lemma \ref{B(n;0,1,0,0)Alperin}, the only possible
$\Ff$-Alperin proper subgroups of $B(3,r;1,0,0)$ are $E_0$ and
$V_0$ if $r$ is odd and $\gamma_2$, $E_0$ and $V_0$ if $r$ is
even. In both cases, according to Lemmas~\ref{B(n)liftings} and~\ref{B(n)restrictions}, if $E_0$ (respectively $V_0$) is
$\Ff$-Alperin, then $\Out_\Ff(E_0)$ contains the diagonal matrix~$-\Id$ (respectively the matrix $\big(\begin{smallmatrix} 1 & 0\\
0 & -1\end{smallmatrix}\big)$) which is not admissible by Lemma
\ref{B(n)liftings}. Therefore we are left with the case of $r=2k$
when $\gamma_2$ is the only $\Ff$-Alperin proper subgroup. But
$\gamma_2\cong \Z/3^{k-1}\times\Z/3^{k-1}$ is normal in
$B(3,r;1,0,0)$ hence the Sylow $3$-subgroup of
$\Out_\Ff(\gamma_2)$ has size $9$, but that contradicts
$\Out_\Ff(\gamma_2)\leq\GL_2(3)$.
\end{proof}

It remains to study the case $B(3,r;0,\gamma,0)$. In this case we
obtain saturated fusion systems with proper $\Ff$-Alperin
subgroups.

\begin{Not}
In what follows we consider the following notation:
\begin{itemize}
\item Fix the following elements in $\Aut(B(3,r;\beta,\gamma,0))$:
\begin{itemize}
\item $\eta$ an element of order two 
which fixes $E_0$  and permutes $E_{1}$ with $E_{-1}$
and such that
projects to $\big(\begin{smallmatrix} 1 & 0 \\ 0 & -1
\end{smallmatrix}\big)$ in $\Out(B(3,r;\beta,\gamma,0))$.
\item $\omega$ an element of order two which commutes with $\eta$,
which projects to $-\Id$ in $\Out(B(3,r;\beta,\gamma,0))$, and such that fixes $E_i$ for $i \in \{-1,0,1\}$.
\end{itemize}
\item By $N\cdot_\gamma W$ we denote an extension of type $N\cdot W$ such that its Sylow $3$-subgroup is isomorphic
to $B(3,r;0,\gamma,0)$.
\end{itemize}
\end{Not}

With all that information now we get the tables with
the possible fusion systems over $B(3,r;0,\gamma,0)$ which are not the normalizer
of the Sylow $3$-subgroup:

\begin{theorem}\label{sfs-B-gamma}
Let $B$ be a rank two $3$-group of maximal nilpotency class of order at least~$3^4$, and let
$(B,\Ff)$ be a saturated fusion system with at least one proper
$\Ff$-Alperin subgroup. Then it must correspond to one of the cases listed in the following tables.
\begin{itemize}
\item If $B\cong B(3,4;0,0,0)$ then the outer automorphism
group of the $\Ff$-Alperin subgroups are in the following table:

\begin{table}[H]
\begin{center}
\begin{tabular}{|c||c|c|c|c|c|c||c|}
\hline
$B$ & $E_0$ & $E_1$ & $E_{-1}$ & $V_0$ & $V_1$ & $V_{-1}$ & $p$-lfg \pb \\
\hline \hline
 & &\phantom{$L_2(3)$} &\phantom{$S_2(3)$} & $\SL_2(3)$ & - & - & $\Ff(3^{4},1)$ \pb \\ \cline{5-8}
 $\langle \omega \rangle$ & - & - & - & - & $\SL_2(3)$ & $\SL_2(3)$ & $\Ff(3^{4},2)$ \pb \\ \cline{5-8}
& & & & $\SL_2(3)$ & $\SL_2(3)$ & $\SL_2(3)$ &  $L_3^{\pm}(q_1)$ \pb
\\ \hline $\langle \eta\omega \rangle$ & $\SL_2(3)$ &
\multicolumn{2}{c|}{-} & - & \multicolumn{2}{c||}{-} & $E_0\cdot_0\SL_2(3)$ \pb \\ \hline
 & & \multicolumn{2}{c|}{ } & - & \multicolumn{2}{c||}{$\SL_2(3)$}  & $\Ff(3^{4},2).2$ \pb \\ \cline{5-8}
 & - & \multicolumn{2}{c|}{-} &  $\GL_2(3)$ & \multicolumn{2}{c||}{-} & $\Ff(3^{4},1).2$ \pb \\ \cline{6-8}
$\langle \eta, \omega \rangle$  &  & \multicolumn{2}{c|}{ } &   & \multicolumn{2}{c||}{$\SL_2(3)$}  & $L_3^{\pm}(q_1):2$\pb \\ \cline{2-8}
  & $\GL_2(3)$ & \multicolumn{2}{c|}{-} & - & \multicolumn{2}{c||}{-} & $E_0\cdot_0\GL_2(3)$ \pb \\ \cline{6-8}
  &  & \multicolumn{2}{c|}{ } & & \multicolumn{2}{c||}{$\SL_2(3)$} & ${}^3D_4(q_2)$\pb \\ \hline
\end{tabular} \\ \mbox{ }
\caption{s.f.s. over $B(3,4;0,0,0)$.\label{ta:B(3,4;0,0,0)}}
\end{center}
\end{table}

Where $q_1$ and $q_2$ are prime powers such that $\nu_3(q_1 \mp 1)=2$ and $\nu_3(q_2^2-1)=1$.

\item If $B\cong B(3,4;0,2,0)$ then the outer automorphism
group of the $\Ff$-Alperin subgroups are in the following table:

\begin{table}[H]
\begin{center}
\begin{tabular}{|c||c|c||c|}
\hline $B$ &  $E_0$ & $V_0$ & $p$-lfg \pb \\ \hline\hline $\langle
\omega \rangle$ & -- & $\SL_2(3)$ & $\Ff(3^4,3)$ \pb\\
\hline $\langle \eta \omega \rangle$ & $\SL_2(3)$ & -- & $E_0 \cdot_2 \SL_2(3)$ \pb\\
\hline $\langle \eta,\omega \rangle$ & -- & $\GL_2(3)$ & $\Ff(3^4,3).2$\pb\\
\cline{2-4} & $\GL_2(3)$ &  -- \pb & $E_0 \cdot_2 \GL_2(3)$\\ \hline
\end{tabular} \\ \mbox{ }
\caption{s.f.s. over $B(3,4;0,2,0)$.\label{ta:B(3,4;0,2,0)}}
\end{center}
\end{table}

\item If $B\cong B(3,2k;0,0,0)$ with $k\geq 3$ then the outer automorphism
group of the $\Ff$-Alperin subgroups are in the following table:

\begin{table}[H]
\begin{center}
\begin{tabular}{|c||c|c|c|c|c|c||c|}
\hline
$B$ & $E_0$ & $E_1$ & $E_{-1}$ & $V_0$ & $V_1$ & $V_{-1}$ & $p$-lfg \pb \\
\hline \hline
 & &\phantom{$E_{-1}$} &\phantom{$S_2$} & $\SL_2(3)$ & - & - & $\Ff(3^{2k},1)$ \pb \\ \cline{5-8}
 $\langle \omega \rangle$ & - & - & - & - & $\SL_2(3)$ & $\SL_2(3)$ & $\Ff(3^{2k},2)$ \pb \\ \cline{5-8}
& & & & $\SL_2(3)$ & $\SL_2(3)$ & $\SL_2(3)$ &  $L_3^{\pm}(q_1)$ \pb
\\ \hline $\langle \eta\omega \rangle$ & $\SL_2(3)$ &
\multicolumn{2}{c|}{-} & - & \multicolumn{2}{c||}{-} & $3\cdot_0\PGL_3(q_2)$ \pb \\ \hline
 & & \multicolumn{2}{c|}{ } & - & \multicolumn{2}{c||}{$\SL_2(3)$}  & $\Ff(3^{2k},2).2$ \pb \\ \cline{5-8}
 & - & \multicolumn{2}{c|}{-} &  $\GL_2(3)$ & \multicolumn{2}{c||}{-} & $\Ff(3^{2k},1).2$ \pb \\ \cline{6-8}
$\langle \eta, \omega \rangle$  &  & \multicolumn{2}{c|}{ } &   & \multicolumn{2}{c||}{$\SL_2(3)$}  & $L_3^{\pm}(q_1):2$\pb \\ \cline{2-8}
  & $\GL_2(3)$ & \multicolumn{2}{c|}{-} & - & \multicolumn{2}{c||}{-} & $3\!\cdot_0\!\PGL_3(q_2)\!\cdot\! 2$ \pb \\ \cline{6-8}
  &  & \multicolumn{2}{c|}{ } & & \multicolumn{2}{c||}{$\SL_2(3)$} & ${}^3D_4(q_3)$\pb \\ \hline
\end{tabular} \\ \mbox{ }
\caption{s.f.s. over $B(3,2k;0,0,0)$ with $k\geq 3$.\label{ta:B(3,2k;0,0,0)}}
\end{center}
\end{table}

Where $q_i$ are prime powers such that $\nu_3(q_1 \mp 1)=k$, $\nu_3(q_2-1)=k-1$ and $\nu_3(q_3^2-1)=k-1$.

\item If $B\cong B(3,2k;0,\gamma,0)$ with $k\geq3$ and $\gamma=1,2$, then the outer automorphism
group of the $\Ff$-Alperin subgroups are in the following table:

\begin{table}[H]
\begin{center}
\begin{tabular}{|c||c|c||c|}
\hline $B$ &  $E_0$ & $V_0$ & $p$-lfg \pb \\ \hline\hline $\langle
\omega \rangle$ & -- & $\SL_2(3)$ & $\Ff(3^{2k},2+\gamma)$ \pb\\
\hline $\langle \eta \omega \rangle$ & $\SL_2(3)$ & -- & $3\cdot_\gamma\PGL_3(q)$\pb\\
\hline $\langle \eta,\omega \rangle$ & -- & $\GL_2(3)$ & $\Ff(3^{2k},2+\gamma).2$\pb\\
\cline{2-4} & $\GL_2(3)$ &  -- \pb & $3\cdot_\gamma\PGL_3(q)\cdot 2$\\ \hline
\end{tabular} \\ \mbox{ }
\caption{s.f.s. over $B(3,2k;0,\gamma,0)$ with $\gamma=1,2$ and $k\geq3$.\label{ta:B(3,2k;0,gamma,0)}}
\end{center}
\end{table}

Where $q$ is a prime power such that $\nu_3(q-1)=k-1$.

\item If $B\cong B(3,2k+1;0,0,0)$ with $k\geq2$ then the outer automorphism
group of the $\Ff$-Alperin subgroups are in the following table:

\begin{table}[H]
\begin{center}
\begin{tabular}{|c||c|c|c|c|c|c|c||c|}
\hline
$B$ & $V_0$ & $V_1$ & $V_{-1}$ & $E_0$ & $E_1$ & $E_{-1}$ & $\gamma_1$ & $p$-lfg \pb \\
\hline \hline
& & \phantom{$V_{-1}$} &  & $\!\SL_2(3)\!$ & - & - & - &  $3.\Ff(3^{2k},1)$\pb \\ \cline{5-9}
 $\langle \omega \rangle$ & - & - & - & - & $\SL_2(3)$ & $\SL_2(3)$ & - & $3.\Ff(3^{2k},2)$\pb \\ \cline{5-9}
& & & & $\!\SL_2(3)\!$ & $\SL_2(3)$ & $\SL_2(3)$ & - & $3\cdot L_3^{\pm}(q_1)$\pb \\ \hline
$\langle \eta \rangle$ & - & \multicolumn{2}{c|}{-} & - & \multicolumn{2}{c|}{-} & $\SL_2(3)$ & $\gamma_1 : \SL_2(3) $\pb \\ \hline
$\!\langle \eta\omega \rangle\!$ & $\!\SL_2(3)\!$ & \multicolumn{2}{c|}{-} & - & \multicolumn{2}{c|}{-} & - & $\PGL_3(q_2)$\pb \\ \hline
 & & \multicolumn{2}{c|}{ } &  & \multicolumn{2}{c|}{-}  & $\GL_2(3)$ & $\gamma_1:\GL_2(3)$ \pb \\ \cline{6-9}
 & & \multicolumn{2}{c|}{ } & - & \multicolumn{2}{c|}{$\SL_2(3)$}  & - & $\!3.\Ff(3^{2k},2).2\!$ \pb \\ \cline{8-9}
 & & \multicolumn{2}{c|}{ } &  & \multicolumn{2}{c|}{}  & $\GL_2(3)$ & $\Ff(3^{2k+1},1)$ \pb \\ \cline{5-9}
 & - & \multicolumn{2}{c|}{-} &   & \multicolumn{2}{c|}{-} & - & $\!3. \Ff(3^{2k},1).2\!$ \pb \\ \cline{8-9}
  & & \multicolumn{2}{c|}{ } & $\!\GL_2(3)\!$ & \multicolumn{2}{c|}{}  & $\GL_2(3)$ & ${}^2F_4(q_3)$\pb \\ \cline{6-9}
$\!\langle \eta, \omega \rangle\!$  
  &  & \multicolumn{2}{c|}{ } &   & \multicolumn{2}{c|}{$\SL_2(3)$}  & - & $\!3\cdot\! L_3^{\pm}(q_1):2\!$ \pb \\ \cline{8-9}
  & & \multicolumn{2}{c|}{ } &  & \multicolumn{2}{c|}{}  & $\GL_2(3)$ & $\Ff(3^{2k+1},2)$\pb \\ \cline{2-9}
  &  & \multicolumn{2}{c|}{ } &  & \multicolumn{2}{c|}{-} & - & $\!\PGL_3(q_2)\cdot\! 2\!$\pb \\ \cline{8-9}
   & $\!\GL_2(3)\!$ & \multicolumn{2}{c|}{-} & -  & \multicolumn{2}{c|}{}  & $\GL_2(3)$ & $\Ff(3^{2k+1},3)$\pb \\ \cline{6-9}
  &  & \multicolumn{2}{c|}{ } & & \multicolumn{2}{c|}{$\SL_2(3)$} & - & $3\cdot{}^3D_4(q_4)$\pb \\ \cline{8-9}
   & & \multicolumn{2}{c|}{ } &  & \multicolumn{2}{c|}{}  & $\GL_2(3)$ & $\Ff(3^{2k+1},4)$\pb \\ \hline
\end{tabular} \\ \mbox{ }
\caption{s.f.s. over $B(3,2k+1;0,0,0)$ with $k\geq 2$.\label{ta:B(3,2k+1;0,0,0)}}
\end{center}
\end{table}
Where $q_i$ are prime powers such that $\nu_3(q_1 \mp 1)=k$,
$\nu_3(q_2 - 1)=k$, $\nu_3(q_3^2-1)=k$ and $\nu_3(q_4^2-1)=k-1$.

\item If $B\cong B(3,2k+1;0,1,0)$ with $k\geq2$ then $\Ff$ is the saturated fusion system associated
to the group $\gamma_1:\SL_2(3)$.
\end{itemize}
The last column of all this tables gives either the information about the groups which have
the corresponding fusion system, either a name encoded as $\Ff(3^r,i)$ to refer to an exotic
$3$-local finite group or the expression of that $3$-local finite group as extension of another
$3$-local finite group.
\end{theorem}

\begin{remark}
As particular cases of the classification we find the exotic fusion systems
over 3-groups of order $3^4$ which were announced previously by
Broto-Levi-Oliver in \cite[Section~5]{BLOs}.
\end{remark}

\begin{proof}
We divide the proof in three parts:

\noindent\textbf{Classification:} In this part we describe the
different possibilities for the saturated fusion system $(B,\Ff)$
by means of the $\Ff$-Alperin subgroups and their outer
automorphisms groups $\Out_\Ff(P)$. Because $\Aut_P(P)\leq
\Aut_\Ff(P)$ by Definition \ref{def:fusion_system},
$\Out_\Ff(P)=\Aut_\Ff(P)/\Aut_P(P)$ also determines $\Aut_\Ff(P)$,
and by Theorem \ref{teo:alperin}, these subgroups of automorphisms
describe completely the category $\Ff$.

By hypothesis we have a proper $\Ff$-Alperin subgroup in
$B(3,r;\beta,\gamma,\delta)$, so by Lemma
\ref{B(n;0,b,c,d)resistant} $\delta=0$, and using now Lemma
\ref{Thm:B(r;0,1,0,0)resistant} also $\beta=0$. So we just have to
cope with $B(3,r;0,\gamma,0)$.

First of all obverve that fixed $i\in\{-1,0,1\}$, $E_i$ and $V_i$
cannot be at the same time $\Ff$-Alperin subgroups: if $E_i$ is
$\Ff$-Alperin then $V_i$ is $\Ff$-conjugate to $\langle
\zeta,\zeta'\rangle$, so $V_i$ is not $\Ff$-centric.

Notice that a saturated fusion systems with
$\Out_\Ff(B(3,r;0,\gamma,0))=1$ cannot contain any proper
$\Ff$-Alperin subgroup: if it had a proper $\Ff$-Alperin subgroup,
by Lemma \ref{B(n)liftings} we would have a nontrivial morphism in
$\Out_\Ff(B(3,r;0,\gamma,0))$.

We now begin the analysis depending on the parity of $r$.

\textbf{Case $r=2k$.} Suppose that $\gamma=0$. If
$\Out_\Ff(B)=\langle \omega \rangle$ then it is immediate from
Lemmas~\ref{B(n)restrictions} and~\ref{B(n)liftings} that $V_i$
may be $\Ff$-Alperin but not $E_i$. The reason is that
$\big(\begin{smallmatrix} -1 & 0\\ 0 & -1\end{smallmatrix}\big)\in
\SL_2(3)$ lifts to $\omega\in\Out_\Ff(B)=\langle \omega \rangle$
from $\Out_\Ff(V_i)$ (it is admissible and the lifting does
exist), while it would lift to $\eta\omega \notin \Out_\Ff(B)$
from $\Out_\Ff(E_0)$ (it is admissible but the lifting does not
exist), and it does not lift for $E_i$ with $i=-1,1$ (it is not
admissible). It is also deduced from Lemmas \ref{B(n)restrictions}
and \ref{B(n)liftings} that no $V_i$ can be $\Ff$-conjugate to
$V_j$ for $i\neq j$ because inner conjugation in $B$ does not move
$\{V_{-1},V_0,V_1\}$ and outer conjugation is induced just by
$\eta$ and $\eta\omega$. If $V_i$ is $\Ff$-Alperin then
$\Out_\Ff(V_i)$ must equal $\SL_2(3)$ because otherwise we would
have a nontrivial element in $\Out_\Ff(B)$ different from $\omega$
for $V_0$, or we would have a non admissible map in $E_{-1}$ or
$E_1$ for $V_{-1}$ or $V_1$ respectively. So one, two or three of
the $V_i$'s can be $\Ff$-Alperin. A symmetry argument, obtained by
conjugation by $B$ (see remarks after Lemma
\ref{B(n)F-automorphisms}), yields the first three rows of tables
\ref{ta:B(3,4;0,0,0)} and \ref{ta:B(3,2k;0,0,0)}. Now assume
$\Out_\Ff(B)=\langle \eta \omega \rangle$. Looking at Lemmas
\ref{B(n)F-automorphisms}, \ref{B(n)restrictions} and
\ref{B(n)liftings} we obtain that only $E_0$ can be $\Ff$-Alperin
and moreover, $\Out_\Ff(E_0)=\SL_2(3)$ because otherwise
$\Out_\Ff(B)$ would be $\Z/2\times\Z/2$. This is the fourth row of
tables \ref{ta:B(3,4;0,0,0)} and \ref{ta:B(3,2k;0,0,0)}. For
$\Out_\Ff(B)=\langle \eta \rangle$ there is no chance for $E_i$ or
$V_i$ to be $\Ff$-Alperin. It remains to cope with the case
$\Out_\Ff(B)=\langle \eta, \omega \rangle$. First, from the
argument above $E_1$ and $E_{-1}$ cannot be $\Ff$-Alperin. Second,
recall from the beginning of this proof that, for $i$ fixed, $E_i$
and $V_i$ cannot be $\Ff$-Alperin simultaneously. Third, notice
that $\eta$ (and $\eta\omega$) swaps $V_1$ and $V_{-1}$, so they
are $\Ff$-conjugate. Lastly, from Lemma \ref{B(n)F-automorphisms},
the only possibility for the outer automorphism groups of $E_0$
and $V_0$ is $\GL_2(3)$ in case they are $\Ff$-Alperin.
Analogously if $V_i$ is $\Ff$-Alperin, for $i=\pm 1$, then
$\Out_\Ff(V_i)$ must be $\SL_2(3)$. Now a case by case checking
yields the last five entries of tables \ref{ta:B(3,4;0,0,0)} and
\ref{ta:B(3,2k;0,0,0)}. If $r>4$ and $\gamma=1,2$ (or $r=4$ and
$\gamma=2$) recall from the conditions in the table of Lemma
\ref{B(n;0,0,c,0)Alperin} that only $E_0$ and~$V_0$ are allowed to
be $\Ff$-Alperin. Similar arguments to those above lead us to
tables \ref{ta:B(3,4;0,2,0)} and \ref{ta:B(3,2k;0,gamma,0)}.

\textbf{Case $r=2k+1$.} For $\gamma=0$ the fusion systems in the
table \ref{ta:B(3,2k+1;0,0,0)} are obtained by similar arguments
to those of the preceding cases, bearing $\gamma_1$ may be
$\Ff$-Alperin too in mind. To fill in this table, notice that if
some $E_i$ or $V_i$ is $\Ff$-Alperin then $-\Id$ is an outer
automorphism of this group that must lift, following Lemma
\ref{B(n)restrictions}, to $\omega$ or $\eta\omega$ for $E_i$ and
$V_i$ respectively. But $\omega$ and $\eta\omega$ restrict in
$\Out_\Ff(\gamma_1)$ to automorphisms of determinant $-1$, so in
case $\gamma_1$ is $\Ff$-Alperin, when some $E_i$ or $V_i$ is so,
its outer automorphism group must be
$\Out_\Ff(\gamma_1)=\GL_2(3)$. Notice that when $\gamma_1$ is
$\Ff$-Alperin, $\Out_\Ff(B)$ must contain $\eta$ by Lemma
\ref{B(n)restrictions} and so $E_{-1}$ and~$E_1$ are
$\Ff$-conjugate. In case of $\gamma=1$, Lemmas
\ref{B(n)F-automorphisms} and \ref{B(n)automorphisms} imply that
only $\big(\begin{smallmatrix} 1 & 0\\ 0 &
-1\end{smallmatrix}\big)$ can be in $\Out_\Ff(B)$. If $\gamma_1$,
$E_0$ or $V_0$ were $\Ff$-Alperin then $\Out_\Ff(B)$ would contain
$\eta$, $\omega$ or $\eta\omega$ respectively. So the only chance
is that $\gamma_1$ is the only $\Ff$-Alperin subgroup. Notice also
that $\Out_\Ff(\gamma_1)$ must equal $\SL_2(3)$, because in other
case there would be a non-trivial element in $\Out_\Ff(B)$
different from~$\eta$.

\noindent\textbf{Saturation:} Now we prove that all the fusion
systems obtained in the classification part of this proof are
saturated by means of \cite[Proposition 5.3]{BLOs}. To show that a
certain fusion system $(B,\Ff)$ appearing in the tables is
saturated the method consists in setting $G\definicio
B:\Out_\Ff(B)$, where $\Out_\Ff(B)$ is the entry for $B$ in the
table, and for each $G$-conjugacy class of $\Ff$-Alperin subgroups
choosing a representative $P\leq B$ and setting
$K_P\definicio\Out_G(P)$ ($K_P$ is determined by Lemma
\ref{B(n)restrictions}) and $\Delta_P\definicio \Out_\Ff(P)$,
where $\Out_\Ff(P)$ is the entry for $P$ in the table. In order to
obtain that the fusion system under consideration is saturated,
for each chosen $\Ff$-Alperin subgroup $P\leq B$ it must be
checked that:
\begin{enumerate}
\item $P$ does not contain any proper $\Ff$-centric subgroup. \item
$p \nmid [\Delta_P:K_P]$ and for each $\alpha\in\Delta_P \setminus
K_P$, $K_P\cap\alpha^{-1}K_P\alpha$ has order prime to $p$.
\label{saturado_condicion_2}
\end{enumerate}
On the one hand, by Lemma \ref{B(n;0,0,c,0)Alperin} the first condition is fulfilled by all the
fusion systems obtained in the classification part of this proof. On the other hand, it is
verified
that $\Out_B(P)$ equals $\langle \big(\begin{smallmatrix} 1 & 0\\
1 & 1\end{smallmatrix}\big)\rangle$, $\langle \big(\begin{smallmatrix} 1 & 1\\
0 & 1\end{smallmatrix}\big)\rangle$ and $\langle \big(\begin{smallmatrix} 1 & 1\\
0 & 1\end{smallmatrix}\big)\rangle$ for $\gamma_1$, $E_i$ and $V_i$ respectively. Denoting by
$\mu_P$ this order~$3$ outer automorphism it is an easy check that for all the fusion systems in
the tables:
\begin{itemize}
\item $P$ is $\Ff$-Alperin with $\Delta_P=\SL_2(3)$ just in case
$K_P=\langle \mu_P, \big(\begin{smallmatrix} -1 & 0\\
0 & -1\end{smallmatrix}\big)\rangle$, \item $P$ is $\Ff$-Alperin
with $\Delta_P=\GL_2(3)$ just in case $K_P=\langle
\mu_P, \big(\begin{smallmatrix} -1 & 0\\
0 & 1\end{smallmatrix}\big),$ $\big(\begin{smallmatrix} 1 & 0\\
0 & -1\end{smallmatrix}\big)\rangle$.
\end{itemize}
Now the two pairs $(K_P,\Delta_P)$ above verify the condition
(\ref{saturado_condicion_2}).\\

Notice that the group $G=B:\Out_\Ff(B)$ defined earlier can be
constructed because $3$ does not divide the order of $\Out_\Ff(B)$
(see Remark \ref{rmk:orderW}) and the projection
$\Aut(B)\twoheadrightarrow\Out(B)$ has kernel a $3$-group, and
thus there is a lifting of $\Out_\Ff(B)$ to $\Aut(B)$. A more
delicate point in the proofs of classification and saturation is
that (recall the remarks after the Lemma
\ref{B(n;0,1,0,0)Alperin}) the outer automorphisms groups
$\Out_\Ff(P)$ for $P\leq B$ $\Ff$-Alperin ($P$ can be the whole
$B$) appearing in the tables are described as subgroups of
$\GL_2(3)$, and that while for $P=E_i,V_i$ the Frattini maps
$\Out(P)\stackrel{\rho}\rightarrow \GL_2(3)$ are isomorphisms, for
$\gamma_1$ and $B$ they are not.

The choice of $\SL_2(3)$ and $\GL_2(3)$ lying in
$\Out(\gamma_1)=\Aut(\gamma_1)$ and of $\Z/2$ and $\Z/2\times
\Z/2\leq \Out(B)$ are not totally arbitrary. In fact, the choice
of $\Aut_\Ff(B)$ must go by the restriction map
$\Aut(B)\rightarrow\Aut(\gamma_1)$ to the choice of
$\Aut_\Ff(\gamma_1)$. Moreover, as $\gamma_1$ is characteristic in
$B$ and diagonal automorphisms in
$\Aut_\Ff(\gamma_1)=\Out_\Ff(\gamma_1)$ must lift to $\Aut_\Ff(B)$
(Lemma \ref{B(n)liftings}), one can check that the choice of
$\Aut_\Ff(\gamma_1)$ determines completely the choice of
$\Out_\Ff(B)$.

Now we prove that different choices gives isomorphic saturated
fusions systems. Let $(B,\Ff)$ and $(B,\Ff')$ be saturated fusion
systems that correspond to the same row in some table of the
classification. Suppose first that $\gamma_1$ is not
$\Ff$-Alperin: the two semidirect products $H=B:\Out_\Ff(B)$ and
$H'=B:\Out_{\Ff'}(B)$ are isomorphic because the lifts of
$\Out_\Ff(B)$ and $\Out_{\Ff'}(B)$ to $\Aut(B)$ are
$\Aut(B)$-conjugate as the projection
$\Aut(B)\twoheadrightarrow\GL_2(3)$ has kernel a $3$-group. Then
we have an isomorphism of categories $\Ff_B(H)\cong\Ff_B(H')$
which can be extended to an isomorphism $\Ff\cong\Ff'$.

If $\gamma_1$ is $\Ff$-Alperin then (recall that
$\Aut_\Ff(\gamma_1)$ determines $\Out_\Ff(B)$) we build the
semidirect products $H=\gamma_1:\Out_\Ff(\gamma_1)$ and
$H'=\gamma_1:\Out_{\Ff'}(\gamma_1)$. These groups are isomorphic
by Corollary \ref{gamma1_semidirectos} and have Sylow $3$-subgroup
$B$. Then we have an isomorphism of categories
$\Ff_B(H)\cong\Ff_B(H')$ which can be extended to an isomorphism
$\Ff\cong\Ff'$.

\noindent\textbf{Exoticism:} To justify the values in the last
column we have to cope with the possible finite groups with the
fusion systems described there.

Consider first all the fusion systems in the tables such that they
have at least one $\Ff$-Alperin rank two elementary abelian
subgroup (respectively at least one $\Ff$-Alperin subgroup
isomorphic to $3^{1+2}_+$ and also $\gamma_1$ is $\Ff$-Alperin).
Consider $N\lneq B(3,r;0,\gamma,0)$ a nontrivial proper normal
subgroup which is strongly closed in $\Ff$. By Lemma
\ref{lem:subgruposnormales}, $N$ must contain the center of $B$,
and as there is an $\Ff$-Alperin rank two elementary subgroup $N$
must also contain~$s$ (respectively, if $\gamma_1$ is
$\Ff$-radical $N$ must also contain $\gamma_{r-2}$, and as there
is an $\Ff$-Alperin subgroup isomorphic to $3^{1+2}_+$, $N$ must
contain $s$ too). Again by Lemma \ref{lem:subgruposnormales} $N$
must be isomorphic to $3^{1+2}_+$ if $r=4$ or $B(3,r-1;0,0,0)$ if
$r>4$.

In all these cases we can apply Proposition
\ref{proposition9.2deBLO2}, getting that if they are the fusion
system of a group $G$, then $G$ can be choosen to be almost
simple. Moreover the $3$-rank of $G$ and the simple group of which
$G$ is an extension must be two, so we have to look at the list of
all the simple groups of $3$-rank two:

\begin{enumerate}[a)]
\item The information about the sporadic simple groups can be deduced from
\cite[Tables 5.3 \& 5.6.1]{GLS}, getting that all of the groups in
that family which $3$-rank equals two have Sylow $3$-subgroup of
order at most $3^3$, and there are not outer automorphisms of
order $3$.

\item The $p$-rank over the Lie type simple groups in a field
of characteristic $p$ are in \cite[Table~3.3.1]{GLS}, where
taking $p=3$ and the possibilities of the groups of $3$-rank two one
gets that the order of the Sylow $3$-subgroup is at most $3^3$, and again
there are not outer automorphism of order $3$.

\item The Lie type simple groups in characteristic prime to $3$
have a unique elementary $3$-subgroup of maximal rank, out of
$L_3(q)$ with $3|q-1$, $L_3^-(q)$ with $3|q+1$, $G_2(q)$,
${}^3D_4(q)$ or ${}^2F_4(q)$ by \cite[10-2]{gorenstein-lyons}. So,
as $3^{1+2}_+$ and $B(3,r;0,\gamma,0)$ do not have a unique
elementary abelian $3$-subgroup of maximal rank, we have to look
at the fusion systems of this small list.

The fusion systems induced by $L_3^+(q)$, when $3|(q-1)$, and by
$L_3^-(q)$, when $3|(q+1)$, are the same, and can be deduced from
\cite[Example 3.6]{BM} and \cite{alperin-fong} using that there is
a bijection between radical subgroups in $\SL_3(q)$ and radical
subgroups in $\PSL_3(q)$, so obtaining the result in Table
\ref{ta:B(3,2k;0,0,0)}. The extensions of $L^{\pm}_3(q)$ by a
group of order prime to $3$ must be also considered, getting
fusion systems over the same $3$-group. Finally,
$\Out(L_3^{\pm}(q))$ has $3$-torsion, so we must consider the
possible extensions, getting the group $\PGL_3(q)$ and an
extension $\PGL_3(q).2$. The study of the proper radical subgroups
in this case is done in~\cite{alperin-fong}, getting that the only
proper $\Ff$-radical is $V_0$.

The fusion system of $G_2(q)$ is studied in \cite{kl2} and \cite{cooperstein}, getting
that it corresponds to the fusion system labeled as $3L_3^+(q):2$.

The fusion system of ${}^3D_4(q)$ can be deduced from \cite{kl}, getting the
desired result.

Finally the fusion system of ${}^2F_4(q)$ has been studied in \cite[Example 9.7]{BM}.
\end{enumerate}

This classification tells us that all the other cases where there is
an $\Ff$-Alperin rank two elementary abelian $3$-subgroup,
and also the ones such that $\gamma_1$ and one subgroup isomorphic to
$3^{1+2}_+$ are $\Ff$-radical,
must correspond to exotic $p$-local finite groups.

Consider now the cases where the only proper $\Ff$-Alperin subgroup is $\gamma_1$. In those cases
it is straightforward to check that they correspond to the groups $\gamma_1 : \SL_2(3)$ and $\gamma_1 :
\GL_2(3)$, where the actions are described in Lemma \ref{lem:acciones}.

In all of the remaining cases, the ones where all the proper
$\Ff$-Alperin subgroups are isomorphic to $3^{1+2}_+$, the
normalizer of the center of 
$B(3,r;0,\gamma,0)$ in $\Ff$ is the
whole fusion system $\Ff$ (i.e. $Z(B(3,r;0,\gamma,0))$ is normal
in $\Ff$).

Consider first the ones where $Z(B(3,r;0,\gamma,0))$ is central in $\Ff$,
i.e. the ones where for all $E_i$ proper $\Ff$-radical
subgroup isomorphic to $3^{1+2}_+$ we have $\Out_\Ff(E_i)=\SL_2(3)$. Using Lemma \ref{lem:extensionescentrales}, we
get that they correspond to groups if and only if the quotient by the
center corresponds also to a group, getting again the results in the tables.

Finally it remains to justify that $3\Ff(3^{2k},1).2$ and $3\Ff(3^{2k},2).2$
are not the fusion system of finite groups. Consider $G$ a finite group with
one of those fusion systems, and consider $Z(B)$ the center of $B(3,2k+1;0,0,0)$.
Consider now the fusion system constructed as the
centralizer of $Z(B)$ (Definition \ref{def:centralizer}) then, by
Remark \ref{rmk:centralizadorygrupo} $3\Ff(3^{2k},1)$ or $3\Ff(3^{2k},2)$
would be also the fusion system of the group $C_{Z(B)}(G)$, and we know that these
are exotic.
\end{proof}

%%%%%%%%%%%%%%%%%%%%%%%%%%%%%%%%%%%%%%%%%%%%%%%%%%%%%%%%%%%%%%%%
% Appendix: group theory
%%%%%%%%%%%%%%%%%%%%%%%%%%%%%%%%%%%%%%%%%%%%%%%%%%%%%%%%%%%%%%%%

\appendix

\section{Rank two $p$-groups}\label{rank-two-p-groups}
In this appendix we recall all the information and properties
of $p$-rank two $p$-groups that
we need to classify the saturated fusion systems over these groups.

The classification of the rank two $p$-groups, $p>2$, traces back
to Blackburn (e.g. see \cite[Theorem 3.1]{leary}):

\begin{theorem}\label{thm:clasificacionrango2}
Let $p$ be an odd prime. Then the $p$-groups of $p$-rank two are
the ones listed here:
\begin{enumerate}[\rm (i)]
\item The non-cyclic metacyclic $p$-groups, which we denote $M(p,r)$.
\item The groups $C(p,r)$, $r \geq 3$ defined by the following presentation:
$$ 
C(p,r) \!\definicio\! \langle a,b,c \mid a^p=b^p=c^{p^{r-2}}=1,
[a,b]=c^{p^{r-3}}, c \in Z(C(p,r)) \rangle . $$
\item The groups $G(p,r;\epsilon)$, where $r \geq 4$ and $\epsilon$ is either $1$ or a
quadratic non-residue modulo~$p$ defined by the following presentation:
$$
G(p,r;\epsilon)\!\definicio\! \langle  a,b,c \mid
a^p\!=\!b^p\!=\!c^{p^{r-2}}\!\!\!\!=\![b,c]\!=\!1,[a,b^{-1}]\!=\!c^{\epsilon p^{r-3}}\!\!,[a,c]\!=\!b
\rangle .
$$
\item If $p=3$ the $3$-groups of maximal nilpotency class, except the
cyclic groups and the wreath product of $\Z/3$ by itself.
\end{enumerate}
Where $[x,y]=x^{-1}y^{-1}x y$.
\end{theorem}
\begin{proof}
According to \cite[Theorem 5.4.15]{gorenstein-68}, the class of
rank two $p$-groups ($p$ odd) agrees with the class of $p$-groups
in which every maximal normal abelian subgroup has rank two, or
equivalently every maximal normal elementary abelian subgroup has
rank two. As the only group of order $p^3$ which requires at least
three generators
is $(\Z/p)^3$, the class of rank two $p$-groups ($p$ odd) agrees
with the class of $p$-groups in which every normal subgroup of
size~$p^3$ is generated by at most two elements. The latter class
of $p$-groups is described in \cite[Theorem~4.1]{BB1} when the
order of the group is $p^n$ for $n\geq 5$ (and $p>2$) while the
case $n\leq 4$ can be deduced for the classification of $p$-groups
of size at most $p^4$ \cite[p.\ 145--146]{burnside}.
\end{proof}

To complete the classification above we also need a description of
maximal nilpotency class $3$-groups, which is given in \cite[last
paragraph p.\ 88]{BB}:

\begin{theorem}\label{thm:maxclass}
The non cyclic $3$-groups of maximal nilpotency class and order greater than~$3^3$
are the groups $B(3,r;\beta,\gamma,\delta)$ with
$(\beta,\gamma,\delta)$ taking the values:
\begin{itemize}
\item For any $r\geq 5$, $(\beta,\gamma,\delta)=(1,0,\delta)$, with
$\delta\in\{0,1,2\}$.
\item For even $r\geq 4$,
$(\beta,\gamma,\delta)\in\{(0,\gamma,0), (0,0,\delta)\}$, with
$\gamma\in\{1,2\}$ and $\delta\in\{0,1\}$.
\item For odd $r\geq 5$,
$(\beta,\gamma,\delta)\in\{(0,1,0), (0,0,\delta)\}$, with
$\delta\in\{0,1\}$.
\end{itemize}
With these parameters, $B(3,r;\beta,\gamma,\delta)$ is the group
of order $3^r$ defined by the set of
generators $\{ s, s_1, s_2, \dots, s_{r-1}\}$ and relations
\begin{align}
 & s_i=[s_{i-1},s] \, \mbox{for $i\in\{2,3, \dots, r-1\}$,} \label{eq:32} \\
 & [s_1, s_2]=s_{r-1}^\beta \, , \label{eq:33} \\
 & [s_1,s_i]=1 \, \mbox{for $i \in \{3,4, \dots, r-1\}$,} \label{eq:35} \\
 & s^3=s_{r-1}^\delta \, ,  \label{eq:36} \\
 & s_1^3 s_2^3 s_3 = s_{r-1}^\gamma \, ,  \label{eq:37}  \\
 & s_i^3 s_{i+1}^3 s_{i+2}=1 \, \mbox{for $i \in \{2,3,\dots,r-1\}$, and assuming $s_r=s_{r+1}=1$}. \label{eq:38}
\end{align}
\end{theorem}

\begin{remark}\label{rmk:idBiG}
For $p=3$ and $r=4$ we have that $B(3,4;0,0,0)\cong G(3,4;1)$, $B(3,4;0,2,0) \cong G(3,4;-1)$ and
$B(3,4;0,1,0)$ is the wreath product $3\wr3$ that has $3$-rank three.
\end{remark}

\begin{remark}
In \cite{BB} the classification of the $p$-groups of maximal rank depends
on four parameters $\alpha$, $\beta$, $\gamma$ and $\delta$, but for $p=3$
we have $\alpha=0$. In all the paper, with the notation $B(3,r;\beta,\gamma,\delta)$
we assume that the parameters $(\beta,\gamma,\delta)$ correspond to
the stated in Theorem~\ref{thm:maxclass} for rank two $3$-groups.
\end{remark}

What follows is a description of the group theoretical properties
of the groups listed in Theorems \ref{thm:clasificacionrango2} and
\ref{thm:maxclass}, that are used along the paper.

We begin with the family  $C(p,r)$:

\begin{lemma} \label{lem:dpC(r)}
Consider $C(p,r)$ as in Theorem \ref{thm:clasificacionrango2},
with the same notation for the generators:
\begin{enumerate}[\rm (a)]
\item The center is $\langle c \rangle \cong
\Z/p^{r-2}$.
\item The commutators are determined by
$[a^ib^j,a^sb^t]=c^{(it-sj)p^{r-3}}$.
\item $C(p,r)=Z(C(p,r))\Omega_1(C(p,r))$ and therefore $C(p,r)$ is
isomorphic to the central product $\Z/p^{r-2}\circ C(p,3)$.
\item The restriction of the elements in $\Aut(C(p,r))$ to
$Z(C(p,r))$ and $\Omega_1(C(p,r))$ provides an isomorphism
$$
\Aut(C(p,r)) \cong
\big(\Z/p^{r-3}\times\ASL_2(p)\big) : (p-1) <\Aut(\Z/p^{r-2})\times\Aut(C(p,3))
$$
that gives rise to a group epimorphism $\rho \colon \Out(C(p,r))
\to \GL_2(p)$ which maps a morphism of type $a \mapsto
a^{i}b^{j}c^{k}$, $b \mapsto a^{i'}b^{j'}c^{k'}$ to the matrix
$\big(\begin{smallmatrix} i & i'
\\ j & j'\end{smallmatrix}\big)$.
\item $\Aut(C(p,r))=\Inn(C(p,r))\rtimes\Out(C(p,r))
\cong(\Z/p\times\Z/p)\rtimes(\Z/p^{r-3}\times\GL_2(p))$.
\item For $G=\GL_2(p)$ or $\SL_2(p)$, the group $\Out(C(p,r))$ contains
just one subgroup isomorphic to~$G$ up to conjugation, but in the
case $G=\SL_2(3)$ and $r>3$. The group $\Out(C(3,r))$, $r>3$,
contains two different conjugacy classes of subgroups isomorphic
to $\SL_2(3)$.
\end{enumerate}
\end{lemma}
\begin{proof}
The statement (a) follows from the presentation of $C(p,r)$ while
(b) can be read from \cite[Lemma 1.1]{dietz-priddy}. The central
product given in (c) is obtained by identifying $\Z/p^{r-2}$ with
$\langle c\rangle=Z(C(p,r))$ and $C(p,3)$ with $\langle
a,b\rangle=\Omega_1(C(p,r))$, so their intersection is $\langle
c^{p^{r-3}}\rangle=Z(\langle a,b\rangle)$. Moreover, $Z(C(p,r))$
and $\Omega_1(C(p,r))$ are characteristic subgroups of $C(p,r)$,
and therefore every element in $\Aut(C(p,r))$ maps each of these
subgroups to itself. Then (c) implies that
$$\Aut(C(p,r))=\{(f,g)\in\Aut(\langle c\rangle)\times\Aut(\langle
a,b\rangle)|\quad f(c^{p^{r-3}})=g(c^{p^{r-3}})\}$$ providing the
description of $\Aut(C(p,r))$ in (d) while the morphism $\rho$ is
obtained by considering outer automorphisms and projecting on the
$C(p,3)$ factor (see \cite[Lemma 3.1]{RV}). As
$\Aut(G)=\Inn(G)\rtimes\Out(G)$ for $G=\Z/p^{r-2},C(p,3)$, (d)
implies (e). Finally, notice that if $G=\GL_2(p)$ or $\SL_2(p)$,
$\GL_2(p)$ contains just one copy of $G$ up to conjugacy, and
therefore, the number of subgroups of
$\Out(C(p,r))\cong\Z/p^{r-3}\times\GL_2(p)$ isomorphic to $G$ up
to conjugation depends on $\hom(G,\Z/p^{r-3})$. Since $G$ is
$p$-perfect for $p>3$, the latter set contains just one element
unless $p=3$, $r>3$ and $G=\SL_2(3)$. Finally,
$\hom(\SL_2(3),\Z/3^{r-3})$ contains three elements that give rise
to three subgroups of type $\SL_2(3)$ in $\Out(C(3,r))$,
determined by their Sylow $3$-subgroups $S_i\definicio\langle
(3^{r-4}i,\big(\begin{smallmatrix} 1 & 1 \\ 0 & 1
\end{smallmatrix}\big) )\rangle$ for $i=0,1,2$ (notice that $\SL_2(p)$ is
generated by elements of order $p$ \cite[Theorem
2.8.4]{gorenstein-68}). But $S_1$ and $S_2$ are conjugate in
$\Out(C(3,r))$ and therefore there are just two conjugacy classes
of $\SL_2(3)$ in $\Out(C(3,r))$ if $r>3$.
\end{proof}

\begin{lemma}\label{lem:pcentricr}
Let $H$ be a $p$-centric subgroup of $C(p,r)$, then $H$ is either
the total or $H \cong \Z/p \times \Z/p^{r-2}$.
\end{lemma}
\begin{proof}
Assume $H$ is a $p$-centric centric subgroup of $C(p,r)$. So $H$
must contain the center~$\langle c \rangle$ as a proper subgroup.
Let $a^ib^j \in H\setminus \langle c \rangle$, then we have that
$\langle a^ib^j,c\rangle\cong \Z/p\times\Z/p^{r-2}$ is a
self-centralizing maximal subgroup in $C(p,r)$, so then $H$ is
either the total or $\langle a^ib^j,c\rangle$.
\end{proof}

The following properties of $G(p,r;\epsilon)$ can be deduced directly from
\cite[Section 1]{dietz-priddy}
\begin{lemma} Consider $G(p,r;\epsilon)$ as in Theorem \ref{thm:clasificacionrango2},
with the same notation for the generators:
\begin{enumerate}[\rm (a)]
\item The commutators are determined by the formula
$[a^ib^jc^k,a^sb^tc^u]=b^{iu-sk}c^n$ where
$n=\epsilon p^{r-3}(u\frac{i(i-1)}{2}+js-it-k\frac{s(s-1)}{2})$.
\item The center of $G(p,r;\epsilon)$ is the group generated by $\langle c^p
\rangle$, and, as $r \geq 4$, it contains $\langle c^{p^{r-3}}
\rangle$.
\item There is a unique automorphism in $G(p,r;\epsilon)$ which
maps $\rho(a)=a^ib^jc^{lp^{r-3}}$ and $\rho(c)=b^tc^u$ for any
$i,j,l,t,u \in \{\pm1\}\times \Z/p\times \Z/p\times \Z/p\times
(\Z/p^{r-2})^*$.
\end{enumerate}
\end{lemma}

\begin{lemma}\label{pcentricenG(re)}
If $(p,r;\epsilon) \neq (3,4;1)$, the $p$-centric subgroups of $G(p,r;\epsilon)$
 are the
ones in the following table:
$$
\begin{array}{|c|c|c|}
\hline \pb \mbox{Isomorphism type} & \mbox{Subgroup} \\ \hline
\hline \pb G(p,r;\epsilon) & \langle a,b,c \rangle \\ \pb
\Z/p\times\Z/p^{r-2} & \langle b,c \rangle \\ \pb \Z/p^{r-2} &
\langle ab^ic^j \rangle \mbox{with $i\in\Z/p$ and
$j\in(\Z/p^{r-2})^*$}\\ \pb M(p,r-1) & \langle ac^j,b\rangle
\mbox{with $j\in(\Z/p^{r-2})^*$}\\ \pb \Z/p\times\Z/p^{r-3} &
\langle ab^i,c^p \rangle \\ \pb C(p,r-1) & \langle a,b,c^p \rangle
\\ \hline
\end{array}
$$
\end{lemma}
\begin{proof}
It is clear that the total is a $p$-centric subgroup, so let
$H<G(p,r;\epsilon)$ be a $p$-centric subgroup different from the
total. As it must contain its centralizer in $G(p,r;\epsilon)$ we
have that $\langle c^p \rangle < H$.

We divide the proof in different cases:

\noindent \textbf{Case $H\leq \langle b,c\rangle$:} then, as
$\langle b , c \rangle$
 is commutative,
we have $H=\langle b,c \rangle \cong \Z/p \times \Z/p^{r-2}$, and
using the commutator rules of this group one can check that it is
$p$-centric.

So in the following cases there is an element of the form
$\alpha=a b^i c^j $.

\noindent \textbf{Case $p\nmid j$ and $H$ cyclic:} as $p\nmid j$
we can construct an automorphism of $G(p,r;\epsilon)$ sending $a
\mapsto ab^i$ and $c \mapsto c^j$, so we can compute the order of
$\alpha$ computing the order of $ac$. Now one can check the
following formula: $$ (ac)^n=a^n b^{- {n \choose 2}} c^{n-{n
\choose 3}\epsilon p^{r-3}} \,. $$ So if $(p;r,\epsilon) \neq
(3;4,1)$ we get that $(ac)^p=c^{\pm p}$, so $ac$ has order
$p^{r-2}$ and it is self-centralizing, so it is $p$-centric. In
this case we have $H = \langle \alpha \rangle \cong \Z/p^{r-2}$.

\noindent \textbf{Case $p\nmid j$ and $H$ not cyclic:} we can
assume that $H$ has two generators, and one of them is $\alpha$:
if we consider an element $\beta$ of $H \setminus \langle \alpha
\rangle$ then the order of $\langle \alpha,\beta \rangle$ is at
least $p^{r-1}$, so if we add another element we would have the
total. As we are considering that $H$ is not the total, we can
assume $H=\langle \alpha,\beta\rangle$. We can also assume that
$\beta=b^kc^l$ (if the generator $a$ would appear in the
expression of $\beta$ we could take a power of $\beta$ and
multiply it by $\alpha^{-1}$ to change the generators). So
consider $H=\langle \alpha , \beta \rangle$ with $H$ not cyclic
and different from the total, then we prove that it is metacyclic:
the order of $H$ is $p^{r-1}$ and, as $H$ is not cyclic, then
$\beta \not\in \langle \alpha \rangle$, so either $k\neq 0$ or
$p\nmid l$. If $k=0$ then, as $p \nmid l$ we have $H=\langle
ab^i,c\rangle$, and as $[ab^i,c]=b$, $H$ is the total, which is
not considered. So $k\neq0$ and now we have to distinguish between
two cases: $p | l$ and $p \nmid l$. If $p\nmid l$ we can consider
the inverse of the automorphism $\rho$ in $G(p,r;\epsilon)$ with
$\rho(a)=a$ and $\rho(c)=b^kc^l$ and we have $\rho^{-1}(H)=\langle
ab^sc^t,c\rangle = G(p,r;\epsilon)$, so $H=G(p,r;\epsilon) $, that
implies that $p| l$ and we can consider the group generated by
$\langle ab^ic^j,b c^{pf} \rangle$. One more reduction is
cancellation of $b^i$ by means of a multiplication by
$(bc^{pf})^{-i}$, so the generators are $\alpha=ac^j$ and
$\beta=bc^{pf}$. We can still simplify the generators using that
$\langle\alpha^p\rangle=\langle c^p\rangle$, getting $\langle
ac^j,b\rangle$. To see that it is metacyclic just check that
$[H,H]=\langle c^{p^{r-3}} \rangle \cong \Z/p$, and that the class
of~$\alpha$ in the quotient $H/[H,H]$ has order $p^{r-3}$, the
same order as $H/[H,H]$, so it is cyclic, and $H\cong M(p,r-1)$.

\noindent \textbf{Case $p\mid j$:} we can assume that all the
elements in $H$ are of the form $a^kb^lc^{pm}$, because if there
were an element in $H$ with an exponent in $c$ which is not
multiple of $p$ then we would be in one of the previously studied
cases. As $c^p$ must be in $H$, then we have that $\langle a
b^i,c^p\rangle \leq H$. One can check by means of the commutator
formula that the group $\langle ab^i,c^p\rangle \cong
\Z/p\times\Z/p^{r-3}$ is self-centralizing, so $p$-centric.
Moreover if $H\setminus \langle ab^i,c^p\rangle \neq \emptyset$,
all the other restrictions of this case imply that there is an
element of the form $a^jb^k$ in $H\setminus \langle ab^i,c^p
\rangle$, and an easy calculation gives us that $\langle
ab^i,a^jb^k,c^p\rangle=\langle a, b, c^p \rangle \cong C(p,r-1)$.
This also proves that there is only one $p$-centric subgroup in
$G(p,r;\epsilon)$ isomorphic to $C(p,r-1)$.
\end{proof}

It remains to study maximal nilpotency class $3$-groups, beginning
with the following properties that can be read in \cite{BB} and
\cite[III.\S 14]{huppert}:
\begin{proposition}\label{propos:maxprop}
Consider $B(3,r;\beta,\gamma,\delta)$ as defined in Theorem \ref{thm:maxclass}
with the same notation for the generators.
Then the following hold:
\begin{enumerate}[\rm (a)]
\item From relations {\rm (\ref{eq:32})} to {\rm (\ref{eq:38})} we get:
\begin{equation}\label{eq:40}
(s^{\pm1} s_1^{\zeta_1} \cdots s_{r-1}^{\zeta_{r-1}})^3=
s_{r-1}^{\pm \delta + \gamma \zeta_1 \pm \beta \zeta_1^2} \,.
\end{equation}
\item $\gamma_i(B(3,r;\beta,\gamma,\delta))\definicio \langle s_i,s_{i+1}, \dots
, s_{r-1}\rangle$ are characteristic subgroups of order $3^{r-i}$
generated by $s_i$ and $s_{i+1}
$
 for $i=1,..,r-1$ (assuming $s_{r}=1$).
\item $\gamma_1(B(3,r;\beta,\gamma,\delta))$
is a metacyclic subgroup.
\item $\gamma_1(B(3,r;\beta,\gamma,\delta))$ is abelian if and only if $\beta=0$
.
\item The extension
$$
 {1} \rightarrow \gamma_1 \rightarrow B(3,r;\beta,\gamma,\delta) \stackrel{\pi} \rightarrow
 \Z/3 \rightarrow {1}
$$
is split if and only if $\delta=0$.
\item $Z(B(3,r;\beta,\gamma,\delta))=\gamma_{r-1}(B(3,r;\beta,\gamma,\delta))=\langle s_{r-1} \rangle$.
\end{enumerate}
\end{proposition}

The following lemma is useful when studying the exoticism of
the fusion systems constructed in Section \ref{maximal_class}.

\begin{lemma} \label{lem:subgruposnormales}
Let $N$ be a nontrivial proper normal subgroup in 
$B(3,r;0,\gamma,0)$. Then
\begin{enumerate}[\rm (a)]
\item $N$ contains $Z(B(3,r;0,\gamma,0))$.
\item If $N$  contains $s$ then $N \cong 3^{1+2}_+$ if $r=4$ or $N \cong B(3,r-1;0,0,0)$ if $r>4$.
\end{enumerate}
\end{lemma}
\begin{proof}
According to \cite[Theorem 8.1]{alperin-bell} $N$ must intersect $Z(B(3,r;0,\gamma,0))$
in a nontrivial subgroup.
As in our case the order of the center is $p$, we obtain (a).

From \cite[Lemma 2.2]{BB} we deduce that if the index of $N$ in
$B(3,r;0,\gamma,0)$ is $3^l$ with $l\geq 2$ then $N=\gamma_l$, and
it does not contain $s$. So a proper normal subgroup $N$ which
contains $s$ must be of index $3$ in $B(3,r;0,\gamma,0)$, so
$B(3,r;0,\gamma,0)/N$ is abelian, and the quotient morphism
$B(3,r;0,\gamma,0)\to B(3,r;0,\gamma,0)/N$ factors through
$B(3,r;0,\gamma,0)/\gamma_2(B(3,r;0,\gamma,0))\cong$ 
$\Z/3 \times
\Z/3$ (the commutator of $B(3,r;0,\gamma,0)$ is
$\gamma_2(B(3,r;0,\gamma,0))$), generated by the classes
$\overline{s}$ and~$\overline{s_1}$. Now we have to take the
inverse image of the proper subgroups in $\Z/3\times\Z/3$ of order
$3$, getting that $N$ must be $\gamma_1$, $\langle s,s_2 \rangle$,
$\langle ss_1,s_2 \rangle$ or $\langle ss_1^{-1},s_2 \rangle$.
Only the second contains $s$, getting the second part of the
result.
\end{proof}

In what follows we restrict to $r\geq5$ and the parameters
$(\beta,\gamma,\delta)\in\{(0,\gamma,0),(1,0,0)\}$, with
$\gamma\in\{0,1,2\}$. This includes all the possibilities for the
parameters described in Theorem \ref{thm:maxclass} except the cases
$\delta=1$ and $\{(0,0,0),(0,2,0)\}$ for $r=4$. We fix $k$ the
integer such that $r=2k$ if $r$ is even and $r=2k+1$ if $r$ is odd.

\begin{lemma}\label{lam:ordeness1s2s3}
For the groups $B(3,r;\beta,\gamma,0)$ it holds that:
\begin{enumerate}[\rm (a)]
\item $\gamma_2(B(3,r;\beta,\gamma,0))$ is abelian.
\item The orders of $s_1$, $s_2$ and $s_3$ are $3^k$, $3^k$ and $3^{k-1}$ if $r=2k+1$
and $3^k$, $3^{k-1}$ and $3^{k-1}$ if $r=2k$.
\end{enumerate}
\end{lemma}
\begin{proof}
We check that $\gamma_2(B(3,r;\beta,\gamma,0))$ is abelian. If
$\beta=0$ then $\gamma_1$ is abelian and $\gamma_2 < \gamma_1$,
and if $\beta = 1$ (which implies $\gamma=0$) we have to see that
$s_2$ and $s_3$ commute. We have the following equalities: $$
[s_3,s_2]=s_3^{-1}s_2^{-1}s_3s_2 \stackrel{(\ref{eq:37})}{=}
s_1^3s_2^3s_2^{-1}s_2^{-3}s_1^{-3}s_2=s_1^3s_2^{-1}s_1^{-3}s_2 \,
. $$ So we have reduced to check that $s_2$ and $s_1^3$ commute.

We use equation (\ref{eq:33}) to deduce
$s_2^{-1}s_1s_2=s_1s_{r-1}$, and raise it to the cubic power to
get $s_2^{-1}s_1^3s_2=s_1^3s_{r-1}^3$. Equation~(\ref{eq:38}) for
$i=r-1$ tells us that the order of $s_{r-1}$ is $3$, so $s_2$ and
$s_1^3$ commute.

We now compute the orders of $s_1$, $s_2$ and $s_3$. Begin using
that $s_{r-1}$ has order $3$. Equation~(\ref{eq:38}) for $i=r-2$
yields $s_{r-2}^3s_{r-1}^3=1$, from which $s_{r-2}$ has order $3$
too. For $i=r-3$ the equation becomes
$s_{r-3}^3s_{r-2}^3s_{r-1}=1$, and so $s_{r-3}$ has order $9$. An
induction procedure, taking care of the parity of $r$, provides us
with the desired result.
\end{proof}

The next lemma describes the conjugation by $s$ action on
$\gamma_1$:

\begin{lemma}\label{B(n;0,0,c,0)prop1}
For $B(3,r;\beta,\gamma,0)$ the conjugation by $s$ on the
characteristic subgroup \linebreak
$\gamma_1(B(3,r;\beta,\gamma,0))$ is given by:
$$
\begin{array}{|c|c|c|c|}
\hline \pb  & \mbox{$\beta=0$, $r=2k+1$} & \mbox{$\beta=0$, $r=2k$} & \mbox{$\beta=1$} \\ 
\hline \hline 
\gamma=0 & 
 \begin{array}{c} \pb s_1^s\!=\!s_1s_2 \\  \pb s_2^s\!=\!s_1^{-3}s_2^{-2} \end{array} & 
  \begin{array}{c} \pb s_1^s\!=\!s_1s_2 \\ \pb s_2^s\!=\!s_1^{-3}s_2^{-2} \end{array} & 
   \begin{array}{c} \pb ({{s_1}^f})^s\!=\!{s_1}^f{s_2}^f{s_{r-1}}^{f(1-f)/2} \\ \pb ({{s_2}^f})^s\!=\!{s_1}^{-3f}{s_2}^{-2f} \end{array} \\
\hline   
\gamma=1 &
 \begin{array}{c} \pb s_1^s\!=\!s_1s_2  \\  \pb s_2^s\!=\!s_1^{-3}s_2^{(-3)^{k-1}-2} \end{array} &
  \begin{array}{c} \pb s_1^s\!=\!s_1s_2 \\  \pb s_2^s\!=\!s_1^{-3((-3)^{k-2}+1)}s_2^{-2} \end{array} &
    - \\ 
\hline 
\gamma=2 & 
 - &
  \begin{array}{c} \pb s_1^s\!=\!s_1s_2 \\ \pb s_2^s\!=\!s_1^{3((-3)^{k-2}-1)}s_2^{-2} \end{array} &  
   -  \\
\hline
\end{array}
$$
% where for $B(3,r;0,\gamma,0)$ the matrix $M_s$ in
% the table is with respect to the ordered basis $\{s_1,s_2\}$ and
% we choose the upper sign if $k$ is odd and the lower one if $k$ is even.
\end{lemma}
\begin{proof}
We begin first with the case of $\beta=0$. To find the expression
for the conjugation by $s$ action on $\gamma_1$ notice that
${s_1}^s=s_1[s_1,s]=s_1s_2$ by equation (\ref{eq:32}) and that
analogously ${s_2}^s=s_2[s_2,s]=s_2s_3$. So we need to express
$s_3$ as a product of powers of $s_1$ and $s_2$. We begin writing
$s_{r-1}$ as a product of powers of $s_2$ and $s_3$. Bearing this
objective in mind we use the same equation (\ref{eq:38}) as
before, but beginning with $i=2$. In this case we obtain
$s_4={s_2}^{-3}{s_3}^{-3}$. For $i=3$ the relation is
$$s_5={s_3}^{-3}{s_4}^{-3} ={s_3}^{-3}({s_2}^{-3}{s_3}^{-3})^{-3}
= {s_2}^9{s_3}^6\textit{.}$$ If $i\geq 4$ and we got in an earlier
stage $$s_i={s_2}^{a_i}{s_3}^{b_i}\textit{, }
s_{i+1}={s_2}^{a_{i+1}}{s_3}^{b_{i+1}}$$ then equation
(\ref{eq:38}) reads as 
$$
s_{i+2}\!=\!{s_i}^{-3}{s_{i+1}}^{-3} \!=\!
({s_2}^{a_i}{s_3}^{b_i})^{-3}
({s_2}^{a_{i+1}}{s_3}^{b_{i+1}})^{-3}\!=\!
{s_2}^{-3a_i-3a_{i+1}}{s_3}^{-3b_i-3b_{i+1}}\textit{.}
$$ So
$s_{r-1}={s_2}^{a_{r-1}}{s_3}^{b_{r-1}}$, where $a_{r-1}$ and
$b_{r-1}$ are obtained from the recursive sequences
$a_4=-3$, $a_5=9$, $a_{i+2}=-3a_i-3a_{i+1}$ and
$b_4=-3$, $b_5=6$, $b_{i+2}=-3b_i-3b_{i+1}$ for $i\geq 4$.

Substituting this last result in equation (\ref{eq:37})
in Theorem \ref{thm:maxclass} we reach
$s_3={s_1}^{3b'}{s_2}^{(3-\gamma a_{r-1})b'}$, where $b'(\gamma
b_{r-1}-1)=1$ modulus the order of $s_3$ and
${s_2}^s={s_1}^{3b'}{s_2}^{1+(3-\gamma a_{r-1})b'}$. Finally, some
further calculus using  the recursive sequences $\{a_i\}$ and
$\{b_i\}$ and taking into account separately the three possible
values of $\gamma$ finishes the proof for $\beta=0$.

% In order to obtain the center of $B(3,r;0,\gamma,0)$, an easy
% calculation shows that it equals $\ker(M_s-\Id)$. Now the
% description of $\zeta$ is achieved from the earlier expressions
% for the matrix $M_s$.

For the case $\beta=1$ use the relations in the presentation of
Theorem \ref{thm:maxclass} to find the commutator rules in
$B(3,r;1,0,0)$.
\end{proof}

\begin{remark}\label{rmk:conjs}
In the cases with $\beta=0$ we have that
$\gamma_1(B(3,r;0,\gamma,0))$ is a rank two abelian subgroup
generated by $\{s_1, s_2\}$ so we can identify the conjugations by
$s$ with a matrix $M^{r,\gamma}_s$, obtaining: $$
M^{r,0}_s=\big(\begin{smallmatrix} 1 & -3 \\ 1 &
-2\end{smallmatrix}\big), \quad
M^{2k+1,1}_s=\big(\begin{smallmatrix} 1 & -3 \\ 1 &
(-3)^{k-1}-2\end{smallmatrix}\big), $$ $$
M^{2k,1}_s=\big(\begin{smallmatrix} 1 & -3((-3)^{k-2}+1) \\ 1 &
-2\end{smallmatrix}\big) \mbox{ and }
M^{2k,2}_s=\big(\begin{smallmatrix} 1 & 3((-3)^{k-2}-1) \\ 1 &
-2\end{smallmatrix}\big) \,. $$
\end{remark}

Next we determine the automorphism groups of $B(3,r;0,\gamma,0)$ and $B(3,r;1,0,
0)$.

\begin{lemma}\label{B(n)automorphisms}
The automorphism group of $B(3,r;\beta,\gamma,0)$ consists of the
homomorphisms that send
$$
s \mapsto s^e {s_1}^{e'} {s_2}^{e''}
$$
$$
s_1\mapsto s_1^{f'}  s_2^{f''}
$$ where the parameters verify the following conditions:
$$
\begin{array}{|c|c|c|c|}
\hline \pb & \mbox{$\beta=0$, $r$ odd} & \mbox{$\beta=0$, $r$ even} & \mbox{ $\beta=1$}\\
\hline \hline \pb \gamma=0 & e=\pm1,3\nmid {f'} & e=\pm1,3\nmid {f'} & e=1,3\mid e',3\nmid f'\\
\hline \pb \gamma=1 & e=1,3\mid e',3\nmid {f'}  & e=\pm 1,3\mid e',3\nmid {f'} & - \\
\hline \pb \gamma=2 & - & e=\pm 1,3\mid e',3\nmid {f'} & -\\
\hline
\end{array}
$$
\end{lemma}
\begin{proof}
We begin with the case of $\beta=0$. As $\gamma_1=\langle
s_1,s_2 \rangle$ is characteristic and $B(3,r;0,\gamma,0)$ is
generated by $s$ and $s_1$, every homomorphism $\varphi \in
B(3,r;0,\gamma,0)$ is determined by $$ s \mapsto s^e {s_1}^{e'}
{s_2}^{e''}, s_1\mapsto s_1^{f'}  s_2^{f''} $$ for some integers
$e,e',e'',{f'},{f''}$. In the other way, given such a set of
parameters the equations~(\ref{eq:32})-(\ref{eq:38}) from Theorem
\ref{thm:maxclass} give us which conditions must verify these
parameters in order to get a  homomorphism $\varphi$. In the study
of these equations denote by $M_s$ the matrix of conjugation by
$s$ on $\gamma_1$ described in Remark \ref{rmk:conjs}.
\begin{itemize}
\item Equation (\ref{eq:32}): it is straightforward to check that these
conditions are equivalent to $\varphi(s_i)={s_1}^{{f'_i}}{s_2}^{{f''_i}}$
for every $2\leq i \leq  n-1$ with $\big(\begin{smallmatrix} {f'_i} \\
{f''_i}\end{smallmatrix}\big)= M_e  \big(\begin{smallmatrix} {f'_{i-1}} \\
{f''_{i-1}}\end{smallmatrix}\big)$, where $M_e={M_s}^e-\Id$, ${f'_1}={f'}$
and ${f''_1}={f''}$. In particular we obtain the value of $\varphi(s_2)$,
and so we get the restriction of $\varphi$ to $\gamma_1=\langle
s_1,s_2 \rangle$.
\item Equations (\ref{eq:33}) and (\ref{eq:35}): as
$\beta=0$ $\gamma_1$ is abelian and these equations
are satisfied if so are the conditions for equation (\ref{eq:32}).
\item Equation (\ref{eq:36}): using equation (\ref{eq:40})
from Proposition \ref{propos:maxprop} we find that if $\gamma=0$
there are no additional conditions and if $\gamma=1,2$ it must be
verified that $3\mid e'$.
\item Equation (\ref{eq:37}): the condition is:
$$
(3+3M_e+M_e^2-\gamma {M_e}^{r-2})\big(\begin{smallmatrix} {f'}
\\{f''}\end{smallmatrix}\big) =0.
$$ For $\gamma=0$ it holds that $3+3M_e+M_e^2=0$. So there is not
more conditions on the parameters. If $\gamma=1,2$, using the
integer characteristic polynomial of $M_e$ to compute
${M_e}^{r-2}$ we produce two recursive sequences, similar to the
ones in Lemma \ref{B(n;0,0,c,0)prop1}, which $(r-1)$-th term
equals $(3+3M_e+M_e^2-\gamma {M_e}^{r-2})\big(\begin{smallmatrix}
{f'}
\\{f''}\end{smallmatrix}\big)$. Checking details we obtain just the
condition $3\mid {f'}$ if $r$ is odd, $\gamma=1$ and $e=-1$.
\item
Equation (\ref{eq:38}): it is easy to check that, as $r\geq 5$, the imposed conditions for
$i=2,...,r-1$ are equivalent to:
$$
M_e(3+3M_e+M_e^2)\big(\begin{smallmatrix} {f'}
\\{f''}\end{smallmatrix}\big) =0.
$$ Because $M_s^3=\Id$, in fact it holds that
$M_e(3+3M_e+M_e^2)=0$. So there are not additional conditions on
the parameters.
\end{itemize}
So we have just got the conditions on the parameters so that there
exists a homomorphism~$\varphi$ such that the images on $s$ and
$s_1$ are predetermined by these parameters $e,e',e'',{f'},{f''}$.
It just remains to look out the conditions on $\varphi$ so it is
an automorphism. A quick check shows that the conditions are
$3\nmid e$ and the determinant of $\varphi|\gamma_1$ is invertible
modulus $3$, which is equivalent to $3\nmid {f'}$.

For $B(3,r;1,0,0)$ the arguments are similar using the commutator rules
for powers of $s$, $s_1$ and $s_2$.
\end{proof}
In the next two lemmas we find out which copies of $\Z/3^n \times \Z/3^n$ and $C
(3,n)$
are in $B(3,r;\beta,\gamma,0)$ and how they lie inside this group.
\begin{lemma}\label{B(n;0,0,c,0)centricos}
Let $P$ be a proper $p$-centric subgroup of $B(3,r;0,\gamma,0)$
isomorphic to $\Z/3^n \times \Z/3^n$ or $C(3,n)$ for some $n$.
Then $P$ is determined, up to conjugation, by the following table:
\begin{center}
\begin{tabular}{|@{}c@{}|l|}
\hline
\begin{tabular}{p{23mm}|p{29mm}}
 Isomorphism type & Subgroup (up to conjugation)\\
\end{tabular} & Conditions \\
\hline \hline
\begin{tabular}{p{23mm}|p{29mm}}
\pb $\Z/3^k \times \Z/3^k$ & $\gamma_1=\langle s_1,s_2 \rangle$
\end{tabular}
 & $r=2k+1$.\\
\hline
\begin{tabular}{p{23mm}|p{29mm}}
\pbbb $3^{1+2}_+$ & $E_i \!\definicio\! \langle \zeta,\zeta',s {s_1}^i \rangle$
\\ \hline
\pbbb $\Z/3 \times \Z/3$ & $V_i \definicio \langle \zeta,s {s_1}^i \rangle$
\end{tabular}
&  \parbox{69mm}{$\zeta={s_2}^{3^{k-1}}, \zeta'={s_1}^{3^{k-1}}$ for $r\!=\!2k\!+\!1$, \\
   $\zeta={s_1}^{3^{k-1}}, \zeta'={s_2}^{-3^{k-2}}$ for $r=2k$, \\
   $i\in\{-1,0,1\}$ if $\gamma=0$ and
   $i=0$ if $\gamma=1,2$.}
\\
\hline
\end{tabular}
%\end{tabular}
\end{center}
\end{lemma}
\begin{proof}
Firstly, suppose that $P$ is contained in $\gamma_1$. As $\gamma_1$ is abelian and $P$ must
contain its centralizer as it is $p$-centric, $P$ must equal $\gamma_1$. Recalling the orders of
$s_1$ and $s_2$ we check that only the case $r$ odd is allowed.

Suppose now that $P$ is not contained in $\gamma_1$. Then $P$ fits
in the short exact sequence:
$$
 {1} \rightarrow K \rightarrow P \stackrel{\pi}\rightarrow
 \Z/3 \rightarrow {1},
$$
with $K\leq \gamma_1$. If $K=\Z/3^m$ then, as $P\cong\Z/3^n \times \Z/3^n$
 or
$P\cong C(3,n)$, we get that the only possibility is $m=1$ and
$P\cong\Z/3 \times \Z/3$. Suppose then that $K=\Z/3^m \times
\Z/3^{m'}$. Now, checking cases again, the only chance for $P$ is
to be $C(3,n)$. An easy calculation shows that
$C_{B(3,r;0,\gamma,0)}(P)\cong\Z/3$. As this centralizer must
contain the center of $P\cong C(3,n)$, which from Lemma
\ref{lem:dpC(r)} is $Z(C(3,n))\cong \Z/3^{n-2}$, $n$ must equal
$3$ and $P\cong 3^{1+2}_+$. Now we determine precisely all the
$p$-centric subgroups isomorphic to $3^{1+2}_+$ or $\Z/3 \times
\Z/3$, and their orbits under $B(3,r;0,\gamma,0)$-conjugation.
\begin{itemize}
\item We begin with the \textbf{$p$-centric subgroups isomorphic
to $3^{1+2}_+$}. As they are $p$-centric they must contain the
center of $B(3,r;0,\gamma,0)$ which is equal to $\langle \zeta
\rangle$. So, from the earlier discussion, they are of the form
$3^{1+2}_{+a}\definicio\langle \zeta,\zeta'\rangle : \langle sa
\rangle$, where $\zeta'$ is as in the statement of the lemma and
$a\in \gamma_1$ is, in principle, arbitrary. Now, $sa$ having
order $3$ is seen to be equivalent to $N(a)=0$ where $N$ is the
norm operator $N =1+s+s^2$. On the other hand, from the
description of $M_s$ in Remark \ref{rmk:conjs}, if $a={s_1}^{a_1
}{s_2}^{a_2}$ and $a'={s_1}^{a_1'}{s_2}^{a_2'}$ then
$3^{1+2}_{+a}$ and $3^{1+2}_{+a'}$ are
$B(3,r;0,\gamma,0)$-conjugate if and only if $a_1$ and $a_1'$ are
congruent modulus $3$. Because $N=0$ for $\gamma=0$ and $a_1\equiv
0$ mod $3$ for every $a={s_1}^{a_1}{s_2}^{a_2}\in \ker (N)$ for
$\gamma=1,2$, we obtain the desired results.
\item The argument to obtain all the
\textbf{$p$-centric subgroups isomorphic to
 $\Z/3\times\Z/3$} is similar. First, from the earlier discussion,
 they must be of the form $\langle \zeta,sa \rangle$ for $\zeta$
 as in the statement and some $a\in\gamma_1$. As before, the condition for $\langle sa
\rangle$ to have order $3$ is $N(a)=0$, and $\langle \zeta,sa
\rangle$ and $\langle \zeta,sa' \rangle$ are
$B(3,r;0,\gamma,0)$-conjugate if and only if $a_1$ and~$a_1'$ are
in the same class modulus $3$.
\end{itemize}
\end{proof}

\begin{lemma}\label{B(n;0,1,0,0)centricos}
Let $P$ be a proper $p$-centric subgroup of $B(3,r;1,0,0)$
isomorphic to $\Z/3^n \times \Z/3^n$ or $C(3,n)$ for some $n$.
Then $P$ is determined, up to conjugation, by the following table:
$$
\begin{tabular}{|@{}c@{}|l|}
\hline
\begin{tabular}{p{28mm}|p{26mm}}
Isomorphism type & Subgroup (up to conjugation) \\
\end{tabular} & Conditions\\
\hline \hline
\begin{tabular}{p{28mm}|p{26mm}}
\pb $\Z/3^{k-1}\! \times \Z/3^{k-1}$ & $\gamma_2=\langle s_2,s_3 \rangle$
\end{tabular} & $r=2k$.\\
\hline
\begin{tabular}{p{28mm}|p{26mm}}
\pbb $3^{1+2}_+$ & $E_0\definicio \langle \zeta,\zeta',s \rangle$ \\
\hline
\pbb $\Z/3 \times \Z/3$ & $V_0 \definicio \langle \zeta,s \rangle$ \\
\end{tabular} & \parbox{65mm}{$\!\zeta\!=\!{s_2}^{3^{k-1}}\!\!\!\!$, $\zeta'\!=\!{s_3}^{-3^{k-2}}\!$
 for $r\!=\!2k\!+\!1$,\\
   $\!\zeta\!=\!{s_3}^{3^{k-2}}\!\!$, $\zeta'\!=\!{s_2}^{3^{k-2}}$ for $r\!=\!2k$.}\\
\hline
\end{tabular}
$$
\end{lemma}
\begin{proof}
Consider the following short exact sequence induced by the abelian
characteristic subgroup $\gamma_2=\langle s_2,s_3 \rangle$ of $B(3,r;1,0,0)$:
$$
 {1} \rightarrow \gamma_2 \rightarrow B(3,r;1,0,0) \stackrel{\pi}\rightarrow
 \Z/3\times\Z/3 \cong \langle \overline{s},\overline{s_1} \rangle \rightarrow {
1}.
$$
If $P$ is contained in $\gamma_2$ then, as $\gamma_2$ is abelian and $P$
must contains its centralizer, $P$ must equal $\gamma_2=\langle s_2,s_3\rangle$.

Recalling the orders of $s_2$ and $s_3$ we check that only the case $r$ even is
allowed.

Suppose now that $P$ is not contained in $\gamma_2$ and consider the non-trivial
 subgroup $\pi(P)$:

\noindent \textbf{Case $\pi(P)=\langle \overline{s_1}
\rangle$,$\langle \overline{s}\overline{s_1}\rangle$ or
$\langle\overline{s}^{-1}\overline{s_1}\rangle $:}
 then $P$ fits in a non-split short exact sequence:
$$
 {1} \rightarrow K \rightarrow P \stackrel{\pi}\rightarrow \Z/3\rightarrow {1},
$$
with $K\leq \gamma_2$. Checking cases for $K$ and $P$ we obtain that in any case
this short exact sequence would split, which is a contradiction.

\noindent \textbf{Case $\pi(P)=\langle \overline{s}\rangle$:} then
$P$ is a subgroup of $\tau=\langle s_2,s_3,s \rangle\cong
B(3,r-1;0,0,0)$. Now apply Lemma \ref{B(n;0,0,c,0)centricos} and
notice that conjugation by $s_1$ conjugates the three copies of
$3^{1+2}_+$ and $\Z/3 \times \Z/3$.

\noindent \textbf{Case $\pi(P)=\Z/3\times \Z/3$:} If $P\neq
B(3,r;1,0,0)$ then there exists a maximal proper subgroup $H<B$
containing $P$. As $\gamma_2$ is the Frattini subgroup of
$B(3,r;1,0,0)$, that is, the intersection of the maximal
subgroups, then $\pi(H)=\Z/3$. This is a contradiction with
$\pi(P)=\Z/3\times \Z/3$, and thus $P$ equals $B(3,r;1,0,0)$.
\end{proof}

%%%%%%%%%%%%%%%%%%%%%%%%%%%%%%%%%%%%%%%%%%%%%%%%%%%%%%%%%%%%%%%%%%%%
%%%%SECCIONES
%%%%%%%%%%%%%%%%%%%%%%%%%%%%%%%%%%%%%%%%%%%%%%%%%%%%%%%%%%%%%%%%%%%%
In the proof of Theorem \ref{sfs-B-gamma} we use implicitly some
particular copies of $\SL_2(3)$ and $\GL_2(3)$ lying in
$\Aut(\gamma_1)$. These are characterized by containing a fixed
matrix. In the next lemma we show when they do exist:

\begin{lemma}\label{lem:acciones}
Consider $P\cong\Z/3^k\times\Z/3^k$ and $M_s^{2k+1,\gamma}$ the
matrix defined in Remark \ref{rmk:conjs} for the case
$B(3,2k+1;0,\gamma,0)$. Then:
\begin{itemize}
\item For $\gamma=0$ there is, up to conjugacy, one
copy of $\SL_2(3)$ (respectively $\GL_2(3)$) in  $\Aut(P)$
containing $M_s^{2k+1,0}$.
\item For $\gamma=1$ there is, up to conjugacy, one
copy of $\SL_2(3)$ (respectively none of $\GL_2(3)$) in $\Aut(P)$
containing $M_s^{2k+1,1}$.
\end{itemize}
\end{lemma}
\begin{proof}
As $\SL_2(3)$ and $\GL_2(3)$ are $3$-reduced any copy of these
groups lying in $\Aut(P)$ is a lift of a subgroup of $\GL_2(3)$ by
the Frattini map $\Aut(P)\stackrel{\rho}\rightarrow \GL_2(3)$.

It can be checked that for $k=2$ the statements of the Lemma are
true. Now, call $P'$ to the Frattini subgroup of $P$,
$P'\definicio\Z/3^{k-1}\times\Z/3^{k-1}$, and consider the
restriction map with abelian kernel: $$
1\rightarrow(\Z/3)^4\rightarrow\Aut(P)\stackrel{\pi}\rightarrow\Aut(P')\rightarrow1.
$$ Let $A$ denote $(\Z/3)^4$. We use this exact sequence to prove
the statements by induction on~$k$. We suppose the lemma is true
for $k-1$ and prove it for $k\ge 3$. We take $G=\SL_2(3)$ or
$\GL_2(3)$ and $\gamma=0$ or $1$. Call $L^{k,\gamma}\definicio
M_s^{2k+1,\gamma}$, where $M_s^{2k+1,\gamma}$ is the matrix
defined in Remark~\ref{rmk:conjs}. It is straightforward that
$\pi(L^{k,0})=L^{k-1,0}$ and $\pi(L^{k,1})=L^{k-1,0}$. We prove
the lemma in three steps:

\noindent\textbf{Existence:} We show the existence of the three
stated copies. By hypothesis there exists a lift
$G\stackrel{\sigma}\rightarrow\Aut(P')$ such that $\rho
\sigma=\Id_G$ and with $\mu\definicio
L^{k-1,0}=\big(\begin{smallmatrix} 1 & -3
\\ 1 & -2\end{smallmatrix}\big)\in \sigma(G)$. Write $H\definicio \langle \mu\rangle\in \Syl_3(\sigma(G))$ and take the
pullback twice: $$ \xymatrix{ A \ar@{^{(}->}[r]             &
\Aut(P)        \ar@{>>}[r]^\pi                 &   \Aut(P') \\ A
\ar@{=}[u]  \ar@{^{(}->}[r] &  \pi^{-1}(\sigma(G)) \ar@{^{(}->}[u]
\ar@{>>}[r]^\pi &   \sigma(G) \ar@{^{(}->}[u]\\ A
\ar@{=}[u]\ar@{^{(}->}[r]   &  \pi^{-1}(H)    \ar@{^{(}->}[u]
\ar@{>>}[r]^\pi &   H \ar@{^{(}->}[u] } $$ As
$L^{k,0}=\big(\begin{smallmatrix} 1 & -3
\\ 1 & -2\end{smallmatrix}\big)$ lies in $\pi^{-1}(H)$
the bottom short exact sequence splits and its middle term can be
identified with $A:L^{k,0}$. Notice also that $L^{k,1}$ lies in
$A:L^{k,0}$, there are several lifts of $H$ and the $A$-conjugacy
classes of these lifts are in $1-1$ correspondence with
$H^1(H;A)$. In fact, each section $H\rightarrow\pi^{-1}(H)$,
corresponds, by the earlier identification, to a derivation
$d\colon H\rightarrow A$, such that $d(\mu)=a\mu$ and
$d(\mu^2)=b\mu^2$ (notice that the $\mu$'s inside and outside the~$d$ lie in different automorphism groups).

Recall that we are interested in building lifts
$G\rightarrow\Aut(P)$ containing $L^{k,\gamma}$, for which is
enough to give sections of the middle short exact sequence in the
diagram above which images contains $L^{k,\gamma}$. If this
sequence splits then the $A$-conjugacy classes of its sections are
in $1-1$ correspondence with $H^1(G;A)$ (for clarity we do not
write $H^1(\sigma(G);A)$), and the sections which image contains
$L^{k,\gamma}$ are precisely those which goes by the restriction
map $\res^G_H \colon H^1(G;A)\rightarrow H^1(H;A)$ to the class of
the section induced by $L^{k,\gamma}$ in the bottom short exact
sequence.

In fact, that the middle sequence splits is due to $[G:H]$ being
invertible in $A$ and applying a transfer argument we find that 
$\res^G_H \colon \! H^*(G;A)\rightarrow H^*(H;A)$ is a monomorphism
(and $\cor^H_G \colon \! H^*(H;A)\rightarrow H^*(G;A)$ is an
epimorphism), so the class of the middle sequence goes by
$\res^G_H \colon H^2(G;A)\rightarrow H^2(H;A)$ to the class of the
bottom sequence, which is zero, and must be the zero class too,
that is, the split one.

Finally, the sections $\sigma^0,\sigma^1 \colon
H\rightarrow\pi^{-1}(H)$ which take $\mu$ to $L^{k,0}$ and
$L^{k,1}$ correspond to the identically zero derivation and to
$a=\big(\begin{smallmatrix} 1 & 0
\\ \mp 3^{k-1} & 1\pm 3^{k-1}\end{smallmatrix}\big)$, $b=\big(\begin{smallmatrix} 1 &
0\\ 0 & 1 \mp 3^{k-1}\end{smallmatrix}\big)$ respectively. These
sections are in the image of the restriction map $\res^G_H\colon
H^1(G;A)\rightarrow H^1(H;A)$ if and only if they are in the
$G$-invariants in $H^1(H;A)$. And easy check shows that $z\in
H^1(H;A)$ is a $G$-invariant if and only if $gz=z$ for every $g\in
N_{\sigma(G)}(H)$. Once computed the derivations $\Der(H,A)$ and
the principal derivations $P(H,A)$, we apply the action of $g$ on
$\sigma^0$ and $\sigma^1$ at the cochain level, and check that the
class of $\sigma^0$ is always $G$-invariant and that the class of~$\sigma^1$ is $G$-invariant just for $G=\SL_2(3)$ in
$H^1(H;A)\cong \Der(H,A)/P(H,A)$.

\noindent\textbf{Uniqueness:} Now we show that the three found
copies are unique up to $\Aut(P)$-conjugation, as claimed. We use
induction and the same tools as in the first step. Take two lifts
$\sigma_1,\sigma_2:G\rightarrow \Aut(P)$ containing
$L^{k,\gamma}$. Composing with $\pi$ we obtain two maps from $G$
to $\Aut(P')$ containing $L^{k-1,0}$. They are lifts  if they are
injective, that is, if $A\cap \sigma_i(G)$ is trivial for $i=1,2$.
As $\GL_2(3)$ and $\SL_2(3)$ are $3$-reduced, and as $A$ is a
normal $3$-group, these groups are indeed trivial. So, by the
induction hypothesis, the two lifts arriving at $\Aut(P')$ must be
conjugated by some $g'\in\Aut(P')$ which centralizes $L^{k-1,0}$.

It is a straightforward calculation that the order of the
centralizers $C_{\Aut(P)}(L^{k,\gamma})$ for $\gamma=0,1$ is
$2\cdot 3^{2k-1}$, and that of $A\cap C_{\Aut(P)}(L^{k,\gamma})$
for $\gamma=0,1$ is $9$ (for every $k\geq 2$). Because $\pi$ maps
$C_{\Aut(P)}(L^{k,\gamma})$ to $C_{\Aut(P')}(L^{k-1,0})$, a
element counting argument shows that in fact
$\pi(C_{\Aut(P)}(L^{k,\gamma}))=C_{\Aut(P')}(L^{k-1,0})$.

Thus, there exists $g\in\Aut(P)$ with $\pi(g)=g'$ and such that
$g$ centralizes $L^{k,\gamma}$. Therefore the images of $\sigma_1$
and $\sigma_2'\definicio c_g \circ \sigma_2$ contain
$L^{k,\gamma}$ and have the same image by $\pi$, that is, they
both lie in $A:\sigma_1(G)=\pi^{-1}(\pi \sigma_1(G))$.

Choosing the Sylow $3$-subgroup $H=\langle \mu\rangle$ of $\pi \sigma_1(G)$
we can construct a three rows short exact sequences diagram as
before, and argue using the injectivity of the restriction map
$\res^G_H \colon H^1(G;A) \rightarrow H^1(H;A)$ to obtain the
uniqueness. More precisely, as the two sections $\pi^{-1}\colon\pi
\sigma_1(G)\rightarrow \sigma_1(G)$ and $\pi^{-1}:\pi
\sigma_1(G)\rightarrow \sigma_2'(G)$ of the middle sequence of the
diagram induce the same section in the bottom row, that is, the
one which maps $\mu$ to $L^{k,\gamma}$, they must be in the same
class in $H^1(G;A)$, which means that they are $A$-conjugate.
Therefore $\sigma_1$ and~$\sigma_2$ are $\Aut(P)$-conjugate.

\noindent\textbf{Non existence:} The arguments of the two
preceding parts prove also  the non existence of sections
$\GL_2(3)\rightarrow \Aut(P)$ containing $L^{k,1}$.
\end{proof}

\begin{remark}
A cohomology-free proof of the non existence of copies (lifts) of
$\GL_2(3)$ in $\Aut(\Z/3^k\times\Z/3^k)$ containing
$L^{k,1}=M_s^{2k+1,1}$ runs as follows: if this were the case,
then, as the elements of order $3$ form a single conjugacy class
in $\GL_2(3)$, we would obtain that $L^{k,1}$ and its square are
conjugate, and so would have same determinant and trace. But one
can checks that this is not the case.
\end{remark}

If $H$ and $K$ are two groups and $H$ acts on $K$ by
$\varphi:H\rightarrow \Aut(K)$ then we can construct the
semidirect product $K :_\varphi H$. In fact, if $\psi:H\rightarrow
\Aut(K)$ is another action conjugate to $\varphi$, that is, exists
$\alpha\in \Aut(K)$ such that
$\psi(h)=\alpha^{-1}\circ\varphi(h)\circ\alpha$ for every $h\in
H$, then $K :_\varphi H\cong K :_\psi H$. The lemma above implies:
\begin{corollary}\label{gamma1_semidirectos}
There exists the groups $\gamma_1:\SL_2(3)$ and
$\gamma_1:\GL_2(3)$ where the actions maps
$\big(\begin{smallmatrix} 1 & 0\\0 & -1\end{smallmatrix}\big)$ to
$M_s^{2k+1,\gamma}$, with $\gamma=0,1$ for $\SL_2(3)$ and
$\gamma=0$ for $\GL_2(3)$. Moreover, these semidirect products
with actions as stated are unique up to isomorphism.
\end{corollary}

%%%%%%%%%%%%%%%%%%%%%%%%%%%%%%%%%%%%%%%%%%%%%%%%%%%%%%%%
% Bibliografia
%%%%%%%%%%%%%%%%%%%%%%%%%%%%%%%%%%%%%%%%%%%%%%%%%%%%%%%%

\end{document}